\documentclass{amsart}
\usepackage{amsmath}
  \usepackage{paralist}
  \usepackage{graphics} 
  \usepackage{epsfig} 
 \usepackage[colorlinks=true]{hyperref}
\hypersetup{urlcolor=blue, citecolor=red}

\newtheorem{theorem}{Theorem}[section]
\newtheorem{corollary}[theorem]{Corollary}
\newtheorem{lemma}[theorem]{Lemma}
\newtheorem{Prop}{Proposition}[section]

\newtheorem{remark}[theorem]{Remark}

\newcommand{\sgn}{\mbox{sgn}}
\newcommand{\vct}[1]{{\mathbf #1}}       

\newcommand{\SSS}{{\mathcal{S}}}

\newcommand{\bound}{\partial\Omega}

\newcommand{\abs}[1]{\left\vert {#1} \right\vert}
\newcommand{\norm}[1]{\left\Vert {#1} \right\Vert}

\newcommand{\pd}[2]{\frac{\partial {#1}}{\partial {#2}}}

\newcommand{\Ree}{\mbox{Re\,}}
\newcommand{\dbar}{\bar{\partial}}

\newcommand{\CC}{{\mathcal{C}}}

\newcommand{\F}{{\mathcal{F}}}

\newcommand{\OO}{{\mathcal{O}}}

\newcommand{\A}{{\mathcal{A}}}

\newcommand{\T}{{\mathbf{t}}}
\newcommand{\Om}{\Omega}
\newcommand{\DOm}{\partial\Omega}

\newcommand{\R}{{\mathbb R}}

\newcommand{\C}{{\mathbb C}}

\newcommand{\matr}[1]{\left[\!\begin{array}{l} #1 \end{array}\!\right]}

\newcommand{\kk}[1]{{\overline{#1}}}

\newcommand{\TT}{{\mathcal{T}}}
\newcommand{\MM}{{\mathcal{M}}}

\def\res{\mathop{Res}\limits}
\def\Ree{\mathrm{Re}}
\def\Imm{\mathrm{Im}}

\title[D-bar method for acoustic tomography]{Positive-energy D-bar method for acoustic tomography: a computational study}

\author{M.~V.~de Hoop, M.~Lassas, M.~Santacesaria, S.~Siltanen and J.~P.~Tamminen}

\subjclass{}
 \keywords{}

 \email{mdehoop@rice.edu}
 \email{matti.lassas@helsinki.fi}
 \email{matteo.santacesaria@helsinki.fi}
 \email{samuli.siltanen@helsinki.fi}
\email{janne.tamminen@ttu.ee}
 \address[M.~V.~de Hoop]{Department of Computational and Applied Mathematics (CAAM), Rice University, 6100 Main St. - MS 134, Houston, TX, United States}
 \address[M.~Lassas, M.~Santacesaria and S.~Siltanen]{Department of Mathematics and Statistics, P.O. Box 68 (Gustaf H\"allstr\"omin katu 2b), 00014 University of Helsinki, Finland}
 \address[J.~P.~Tamminen]{Department of Mathematics, Tallinn University of Technology, Ehitajate tee 5, 12616 Tallinn, Estonia}

\begin{document}

\begin{abstract}
A new computational method for reconstructing a potential from the Dirichlet-to-Neumann map at positive energy is developed. The method is based on D-bar techniques and it works in absence of exceptional points -- in particular, if the potential is small enough compared to the energy. Numerical tests reveal exceptional points for perturbed, radial potentials. Reconstructions for several potentials are computed using simulated Dirichlet-to-Neumann maps with and without added noise. The new reconstruction method is shown to work well for energy values between $10^{-5}$ and $5$, smaller values giving better results.
\end{abstract}

\maketitle


%
%
\bigskip

\tableofcontents

\section{Introduction}
\noindent
Indirect measurements can often be accurately modelled using boundary value problems involving partial differential equations (PDE). The need to interpret such measurements leads to inverse problems where one aims to recover spatially varying PDE coefficients from boundary data. Examples include Electrical Impedance Tomography (EIT) and Acoustic Tomography (AT).

The nonlinearity and ill-posedness of the above kind of inverse problems call for specialized solution methods. Ideally, the methodology should provide computationally feasible instructions for reconstructing the coefficient of interest in a noise-robust way, and a regularization analysis for the result. One promising approach is the use of {\em nonlinear Fourier transforms} tailored to the problem at hand. Known as the {\em D-bar method}, it already works for EIT in practice \cite{Novikov1988,Nachman1996,Siltanen2000,Mueller2003,Isaacson2006,Knudsen2009,Mueller2012}, and the theory for AT is well-known \cite{Grinevich1986, Novikov1986b, Novikov1988}. A numerical method for AT was recently presented in \cite{Burov2013, Burov2012} where the classical and generalized scattering amplitude are computed from the measurements as a reconstruction step. In this paper we provide a new numerical implementation of the original D-bar method at positive energy, applicable to AT, where we use the nonlinear Fourier transform computed directly from measurements.


Let $\Omega\subset\R^2$ be the unit disk. Consider a bounded function $\varrho: \Omega\rightarrow (0,\infty)$, which models physical \textit{density} and satisfies $\varrho(x)\geq c_\rho>0$. The \textit{pressure} $p$ satisfies the \textit{reduced acoustic equation} for time-harmonic waves with frequency $\omega$,
\begin{equation}\label{AT}
\nabla\cdot(\frac{1}{\varrho}\nabla p)+\omega^2\kappa p  = 0 \quad \textrm{in }\Om,
\end{equation}
where $\kappa(x)$ is the compressibility and the speed of sound is given by $c = (\kappa\varrho)^{-1/2}$. Given the boundary condition $p=f$ on $\DOm$, the inverse problem of AT is to recover $\varrho$ and $\kappa$ from the Dirichlet-to-Neumann (DN) map
\begin{equation} \label{DN-AT-def}
\Lambda_{\omega,\kappa,\varrho}:\quad f\rightarrow \frac{1}{\varrho}\pd{p}{\nu}\rvert_{\DOm}.
\end{equation}

Assume for simplicity that $\varrho(x)=h$ and $\kappa(x)=k$ for all $x$ near $\partial\Omega$ with some positive constants $h$ and $k$. Define
\begin{equation}
q_0=\varrho^{1/2}\Delta\varrho^{-1/2}-(\omega^2\kappa\varrho-E);
\end{equation}
note that $\mbox{supp}(q_0)\subset\Omega$ since we define the {\em energy} $E$ by
\begin{equation}
  E=k\cdot h\cdot \omega^2>0.
\end{equation}
The substitution $u=\varrho^{-1/2}p$ transforms (\ref{AT}) into the Schr\"odinger equation
\begin{equation}\label{schrodinger-nonzero}
(-\Delta+q)u = 0 \quad \textrm{in }\Omega,
\end{equation}
where  $u=f$ on $\DOm$ and $q=q_0-E$. This boundary-value-problem might be well-posed; then for any $f\in H^{1/2}(\DOm)$ it has a unique weak solution $u \in H^1(\Omega)$. In that case we define the DN map
\begin{equation} \label{DNq-def}
\Lambda_q:\quad H^{1/2}(\DOm)\rightarrow H^{-1/2}(\DOm),\quad f\mapsto \pd{u}{\nu}\rvert_{\DOm},  
\end{equation}
where $\nu$ is the unit outer normal to the boundary. More precisely for $f,g\in H^{1/2}(\DOm)$,
\begin{equation} \label{DNq-def2}
  (\Lambda_q f,g)_{\DOm} = \int_\Om(\nabla u\cdot\nabla v + quv)dx,
\end{equation}
where $u$ is the unique weak solution for the boundary value $f$, and $v\in H^1(\Omega)$ with $v|_{\DOm}=g$.  As mentioned in \cite{Nachman1988} we then have
\begin{equation*}
\Lambda_q = \varrho^{1/2}\Lambda_{\omega,\kappa,\varrho}\varrho^{1/2}-\frac{1}{2}\varrho^{-1}\pd{\varrho}{\nu},
\end{equation*}
so our assumptions imply $\Lambda_q = h\cdot \Lambda_{\omega,\kappa,\varrho}$. Assuming that $h,k$ and $\omega$ are known, the inverse problem of AT then takes the following form: given $\Lambda_q$ and the energy $E$, reconstruct the potential $q_0$. This is called the \textit{Gel'fand-Calder\'on problem} posed by Gel'fand \cite{Gelfand1961} and Calder\'on \cite{Calder'on1980}.

We can solve this problem using the D-bar method based of exponentially behaving Complex Geometric Optics (CGO) solutions first introduced by Faddeev \cite{Faddeev1966} and later rediscovered in the context of inverse boundary-value problems by Sylvester and Uhlmann \cite{Sylvester1987}.  The D-bar method is based on the boundary integral equation proved by R. G. Novikov \cite{Novikov1988}, the D-bar equation discovered by Beals and Coifman \cite{Beals1981}, and the relation of the CGO solution and the potential by R. G. Novikov \cite{Novikov1986b}. See Nachman \cite{Nachman1988a} for a discussion of the D-bar method applied to the AT problem.

EIT and AT are related to the Gel'fand-Calder\'on problem by a transformation resulting to different energies: EIT is a zero-energy problem with $E=0$ and AT is a positive energy problem with $E>0$. In the zero-energy case in 2D, for conductivity-type potentials, Nachman \cite{Nachman1996} proved uniqueness and rigorously justified the D-bar reconstruction. The result was later generalized by Bukhgeim \cite{Bukhgeim2008}, who proved global uniqueness for general potentials at any fixed energy.

The three novelties of this paper are:
\begin{enumerate}
\item We create a numerical algorithm for \textit{Faddeev Green's function} for positive energy $E>0$, a significant extension of the zero-energy case introduced in \cite{Siltanen1999}. This is done in section \ref{GZETA} after which the function will be used throughout the numerical computations.
\item We investigate numerically the \textit{exceptional points} which prevent the straightforward use of the D-bar method for reconstruction. This numerically complements the earlier works \cite{Music2013,Music2014,Grinevich2012,Grinevich2013} focusing on the zero and non-zero energy cases. The results can be found in subsection \ref{sec:exceptional}.
\item We propose a new numerical algorithm for the D-bar method at positive energy and test it to reconstruct potentials using simulated DN-maps with and without added noise. In contrast to other methods, our algorithm is able to do reconstructions at low energies. See subsections \ref{relatedwork} and \ref{sec:comparison} for comparisons with other reconstruction schemes (\cite{Novikov01012013,Burov2012,Bukhgeim2008,lakshtanov2015global}). The algorithm is detailed in section \ref{sec:Dbar-numerics} and tested in subsection \ref{sec:reconstruction}.
\end{enumerate}\bigskip

\section{Preliminaries}

\subsection{Complex Geometric Optics solutions and exceptional points}\label{CGO}
Let $q_0 \in L^{\infty}(\Omega)$ be real-valued, and $E >0$. 

Rewrite \eqref{schrodinger-nonzero} and consider special exponentially growing solutions $\psi(x,\zeta)$ of
\begin{equation}\label{schrode1}
(-\Delta+q_0)\psi(\,\cdot\,,\zeta) = E\psi(\,\cdot\,,\zeta) \quad \mbox{in }\R^2,
\end{equation}
where $q_0$ is extended to the plane by zero, $x = [x_1,x_2]^T \in \R^2$ and $\zeta=[\zeta_1\ \zeta_2]^T\in\C^2$ is a spectral parameter with $\Imm(\zeta)\neq\vct{0}$. The exponential behaviour is then included in the requirement
\begin{equation}\label{exp-behaviour}
e^{-i\zeta\cdot x}\psi(x,\zeta)-1 \to 0, \qquad \text{as } |x| \to +\infty
\end{equation}
where $\zeta^2 = \zeta_1^2+\zeta_2^2 = E$ and $\zeta\cdot x = x_1\zeta_1+x_2\zeta_2$.

The requirement $\zeta^2 = E$ comes from the fact that for large $\abs{x}$ we can write $\exp(i\zeta\cdot x)$ in place of $\psi(x,\zeta)$, and the equation \eqref{schrode1} still has to hold, giving us
\begin{equation}\label{E-requirement}
  -\Delta e^{i\zeta\cdot x} = Ee^{i\zeta\cdot x}\quad\Rightarrow\quad\zeta\cdot\zeta = E.
\end{equation}

If $\zeta\in\R^2$, we have a setting of a physical scattering of a particle with momentum $\zeta$. If $\Imm(\zeta)\neq\vct{0}$, we call the solutions $\psi(x,\zeta)$ \textit{Complex Geometric Optics} or \textit{Faddeev-type} solutions (first introduced in \cite{Faddeev1966}).

In order to precisely define the CGO solutions, we first write
$$
\mu(x,\zeta) = e^{-i\zeta\cdot x}\psi(x,\zeta)
$$
and call it a CGO solution as well. This new function satisfies another differential equation; starting from \eqref{schrode1} we get
\begin{eqnarray*}
\begin{array}{rcl}
(-\Delta+q_0-E)e^{i\zeta\cdot x}\mu(x,\zeta) &=& 0 \nonumber\\
e^{i\zeta\cdot x}(-\Delta-2i\zeta\cdot\nabla+\zeta\cdot\zeta+q_0-E)\mu(x,\zeta) &=& 0 \nonumber\\
             (L_\zeta+q_0)\mu(x,\zeta) &=& 0, 
\end{array}
\end{eqnarray*}
where
\begin{equation*}
L_\zeta := -\Delta-2i\zeta\cdot\nabla.
\end{equation*}
The Green function $g_\zeta(x)$ of the operator $L_\zeta$ is called \textit{Faddeev Green's function} and it is explicitly given by the formula
\begin{equation}\label{green_zeta}
g_\zeta(x) = \frac{1}{4 \pi^2} \int_{\R^2} \frac{e^{iy\cdot x}}{y\cdot y + 2\zeta \cdot y}d y,
\end{equation}
for $x \in \R^2$ and $\Imm(\zeta)\neq\vct{0}$. Thus we define $\mu(x,\zeta)$ as the solution of the following Lippmann-Schwinger type equation
\begin{equation}\label{LS-lambda}
  \mu(x,\zeta)= 1-g_\zeta(x)\ast(q_0(x)\mu(x,\zeta)).
\end{equation}

For a given potential $q_0$, we call $\zeta \in \C^2 \setminus \R^2$ an \textit{exceptional point} if and only if integral equation \eqref{LS-lambda} does not admit a unique solution in $L^{\infty}(\R^2)$.

When $\Imm(\zeta) = \vct{0}$, formula \eqref{green_zeta} and equation \eqref{LS-lambda} make no sense; however, the following limits can be defined:
\begin{gather} \label{defpsig}
\psi_{\gamma}(x,\zeta)=\psi(x,\zeta+i0\gamma), \qquad g_{\zeta}^{\gamma}(x)= g_{\zeta+i0\gamma}(x),\\
\mu_{\gamma}(x,\zeta)=\mu(x,\zeta+i0\gamma),
\end{gather}
where $\zeta,\gamma \in \R^2$, $\zeta^2 = E$, $\gamma^2=1$ and $f(\zeta+i0\gamma) = \lim_{\varepsilon \to 0^+}f(\zeta + i\varepsilon \gamma)$.


Following \cite{Grinevich1986}, we make the change of variables 
\begin{equation}\label{lambda-definition}
z=x_1+ix_2, \quad \lambda = \frac{\zeta_1+i\zeta_2}{\sqrt{E}},\quad \zeta = \matr{(\lambda+\frac{1}{\lambda})\frac{\sqrt{E}}{2} \\ (\frac{1}{\lambda}-\lambda)\frac{i\sqrt{E}}{2}}. 
\end{equation}
We call the new parameter $\lambda$ also the spectral parameter. Depending on whether we use $(x,\zeta)$- or $(z,\lambda)$-notation, in place of $\psi(x,\zeta)$, $\mu(x,\zeta)$ and $g_{\zeta}(x)$ we write $\psi(z,\lambda)$, $\mu(z,\lambda)$ and $g_{\lambda}(z)$ respectively, even if the energy is then omitted. In the numerical part of the paper we clearly indicate which energy level we are using in different numerical tests.

Let $\lambda=r\exp(i\theta)$. Rewriting \eqref{lambda-definition} gives
\begin{equation}\label{lambda-def2}
\zeta = \frac{\sqrt{E}}{2}\left((r+\frac{1}{r})\matr{\cos(\theta)\\\sin(\theta)}+i(r-\frac{1}{r})\matr{\sin(\theta)\\-\cos(\theta)}\right).
\end{equation}
It is easy to see that $r\rightarrow 1$ implies $\Imm(\zeta)\rightarrow\vct{0}$, meaning that the CGO solution goes to the limit of physical scattering. For $|\lambda|=1$ ($r=1$) we then define:
\begin{gather}\label{defmupm}
 \psi_{\pm}(z,\lambda) =\psi(z,\lambda(1\mp 0)), \qquad \mu_{\pm}(z,\lambda) = \mu(z,\lambda(1\mp 0)), \\ \label{defgpm}
g_{\lambda}^{\pm}(z)=g_{\lambda(1\mp 0)}(z),
\end{gather}
where $f(1\mp 0 ) = \lim_{\varepsilon \to 0^+}f(1 \mp \varepsilon)$.


%
%

\subsection{The D-bar equation and the boundary integral equation}\label{sec:dbar}
All of the following is included in the papers \cite{Novikov1992}, \cite{Novikov1988} and the survey \cite{Grinevich2000} with different notation. Recall our assumptions of $q_0$ to be real valued and $E>0$.

Define the differential operators
\begin{equation*}
\partial_w = \frac{1}{2}(\partial_{w_1}-i\partial_{w_2}),\quad \dbar_w = \frac{1}{2}(\partial_{w_1}+i\partial_{w_2}),
\end{equation*}
where $w = w_1+iw_2$, and the exponential functions
\begin{eqnarray}
e_{-\lambda}(z) &:=& \exp\left(-\frac{i\sqrt{E}}{2}(1+\frac{1}{\lambda\kk{\lambda}})(z\kk{\lambda}+\kk{z}\lambda)\right)\label{exp-mlambda},\\
e_{\lambda}(z)  &:=& \exp\left(\frac{i\sqrt{E}}{2}(1+\frac{1}{\lambda\kk{\lambda}})(z\kk{\lambda}+\kk{z}\lambda)\right)\label{exp-lambda}.
\end{eqnarray}
For $\zeta \in \C^2 \setminus \R^2$ with $\zeta^2 = E$, not an exceptional point, we can define, for the corresponding $\lambda$, the non-physical scattering transform by
\begin{equation}\label{scattering-transform}
\vct{t}(\lambda) = \int_{\C}e_{\lambda}(z)q_0(z)\mu(z,\lambda)d\Ree z \, d\Imm z,
\end{equation}
Here $d\Ree z \, d\Imm z$ stands for the standard Lebesgue measure on the plane, i.e. $dx_1dx_2$, since $z = x_1+i x_2$. We have the following symmetry that we will use later:
\begin{equation}\label{tsymmetry}
\vct{t}(1/\kk{\lambda}) = \vct{t}(\lambda).
\end{equation}
Further, we define the functions $h_{\pm}$,
\begin{align} \label{defhpm}
h_{\pm}(\lambda,\lambda',E)&=\left(\frac{1}{2 \pi}\right)^2\int_{\C}\exp\left[-\frac{i}{2}\sqrt{E}(\lambda'\bar z+z/\lambda')\right]\\ \nonumber
&\quad \times q_0(z)\psi_{\pm}(z,\lambda)d\Ree z \, d\Imm z,
\end{align}
for $|\lambda|=|\lambda'|=1$. It is also useful to introduce the following auxiliary functions $h_1, h_2$,
\begin{align} \label{defh1}
h_1(\lambda,\lambda')&=\chi_+\left[-\frac{1}{i}\left(\frac{\lambda'}{\lambda}-\frac{\lambda}{\lambda'}\right) \right]h_+(\lambda,\lambda')\\ \nonumber
&\quad -\chi_+\left[\frac{1}{i}\left(\frac{\lambda'}{\lambda}-\frac{\lambda}{\lambda'}\right) \right]h_-(\lambda,\lambda'),\\ \label{defh2}
h_2(\lambda,\lambda')&=\chi_+\left[-\frac{1}{i}\left(\frac{\lambda'}{\lambda}-\frac{\lambda}{\lambda'}\right) \right]h_-(\lambda,\lambda')\\ \nonumber
&\quad -\chi_+\left[\frac{1}{i}\left(\frac{\lambda'}{\lambda}-\frac{\lambda}{\lambda'}\right) \right]h_+(\lambda,\lambda'),
\end{align}
where $\chi_+$ is the Heaviside step function, and $\rho$, solution of the following integral equations,
\begin{subequations}\label{defrho}
\begin{align} \label{defrho1}
&\rho(\lambda,\lambda')+\pi i\int_{|\lambda''|=1}\rho(\lambda,\lambda'')\chi_+\left[\frac{1}{i}\left(\frac{\lambda'}{\lambda''}-\frac{\lambda''}{\lambda'}\right) \right]\\ \nonumber
&\qquad \times h_1(\lambda'',\lambda')|d\lambda''|=-\pi ih_1(\lambda,\lambda'),\\
&\rho(\lambda,\lambda')+\pi i\int_{|\lambda''|=1}\rho(\lambda,\lambda'')\chi_+\left[-\frac{1}{i}\left(\frac{\lambda'}{\lambda''}-\frac{\lambda''}{\lambda'}\right) \right]\\ \nonumber
&\qquad \times h_2(\lambda'',\lambda')|d\lambda''|=-\pi ih_2(\lambda,\lambda'),
\end{align}
\end{subequations}
for $|\lambda|=|\lambda'|=1$. Here and throughout the paper, $|d\lambda''|$ (or later $|d\lambda'|$) stands for the \textit{surface measure} on $\{\lambda \in \C : |\lambda|=1 \}$.

The function $\mu$ satisfies the following non-local Riemann-Hilbert problem (see \cite{Grinevich1986} and \cite{Novikov1992}). We have:
\begin{equation}\label{D-bar}
\dbar_\lambda\mu(z,\lambda) = \sgn(\abs{\lambda}^2-1)\frac{\vct{t}(\lambda)}{4\pi\kk{\lambda}}e_{-\lambda}(z)\kk{\mu(z,\lambda)},
  \end{equation}
for $\lambda$ not an exceptional point and $|\lambda| \neq 1$,
\begin{align}\label{corrt}
\mu_+(z,\lambda) = \mu_-(z,\lambda)+\int_{|\lambda'|=1}\rho(\lambda,\lambda',z)\mu_-(z,\lambda')|d\lambda'|,
\end{align}
for $|\lambda|=1$, where
\begin{equation}\label{defrhoz}
\rho(\lambda,\lambda',z)= \rho(\lambda,\lambda')\exp\left[\frac{i\sqrt{E}}{2}\left((\lambda'-\lambda)\bar z+\left(\frac{1}{\lambda'}-\frac{1}{\lambda}\right)z\right)\right].
\end{equation}
In addition we have
\begin{gather}
\lim_{|\lambda|\to \infty} \mu(z,\lambda)=1,\\ \label{mu-asymptotics}
\mu(z,\lambda) = 1+\frac{\mu_{-1}^{\infty}(z)}{\lambda}+\OO\left(\frac{1}{|\lambda|}\right),\quad \text{for } |\lambda| \to \infty,\\ \label{q0-recon1}
q_0(z)=2i\sqrt{E}\partial_z\mu_{-1}^{\infty}(z).
\end{gather}

Define the operators
\begin{eqnarray}
\CC:\quad \CC f(z,\lambda) &=& \frac{1}{\pi}\int_{\C}\frac{f(z,w)}{w-\lambda}dw\label{Cauchy},\\
\TT:\quad \TT f(z,\lambda) &=& \sgn(\abs{\lambda}^2-1)\frac{\vct{t}(\lambda)}{4\pi\kk{\lambda}}e_{-\lambda}(z)\kk{f(z,\lambda)}\label{T-def},\\
\MM:\quad  \MM f(z,\lambda) &=& \frac{1}{2\pi i}\int_{\abs{w}=1}\frac{dw}{w-\lambda}\int_{|\lambda'|=1}\!\!\!\! \rho(w,\lambda',z)f_+(z,\lambda')|d\lambda'|,
\end{eqnarray}
where $f_{+}(z,\lambda)$ is the limit of $f(z,\lambda)$ when $\abs{\lambda}\rightarrow 1$ as defined in \eqref{defmupm}. By applying the Cauchy-Green formula to \eqref{D-bar} and \eqref{corrt} the CGO solution $\mu(z,\lambda)$ satisfies the integral equation 
  \begin{equation}
\mu(\cdot,\lambda) = 1-(\CC\TT-\MM)\mu(\cdot,\lambda)\label{IEIS-compact},
  \end{equation}
where for each fixed $z'\in\Omega$ we can solve $\mu(z',\lambda)$.

We now review the reconstruction scheme to obtain the scattering data $\vct{t}(\lambda)$ and $h_{\pm}(\lambda,\lambda')$ from the Dirichlet-to-Neumann data. 

Define the operators
\begin{equation}\label{Slambda-def}
(\SSS_\lambda \phi)(z) := \int_{\DOm}G_\lambda(z-y)\phi(y)ds(y),\quad G_\lambda(z) = e^{\frac{i\sqrt{E}}{2}(\lambda \kk{z}+\frac{z}{\lambda})}g_\lambda(z),
\end{equation}
for $z \in \partial \Omega$, $|\lambda| \neq 1$ and
\begin{equation}\label{Slambdapm-def}
(\SSS_\lambda^{\pm} \phi)(z) := \int_{\DOm}G_\lambda^{\pm}(z-y)\phi(y)ds(y),\quad G^{\pm}_\lambda(z) = e^{\frac{i\sqrt{E}}{2}(\lambda \kk{z}+\frac{z}{\lambda})}g_{\lambda}^{\pm}(z),
\end{equation}
for $z \in \partial \Omega$, $|\lambda|=1$.
The CGO solutions $\psi(z,\lambda)$ and $\psi_{\pm}(z,\lambda)$ satisfy the boundary integral equations \cite{Novikov1988}
\begin{align}\label{BIE-lambda}
&(I+\SSS_\lambda(\Lambda_q-\Lambda_{-E}))\psi(\cdot,\lambda)|_{\DOm} = e^{\frac{i\sqrt{E}}{2}(\lambda\kk{z}+\frac{1}{\lambda}z)}|_{\DOm}, \quad \text{for } |\lambda| \neq 1 \textrm{ (non-exceptional)},\\ \label{BIE-lambdapm}
&(I+\SSS_\lambda^{\pm}(\Lambda_q-\Lambda_{-E}))\psi_{\pm}(\cdot,\lambda)|_{\DOm} = e^{\frac{i\sqrt{E}}{2}(\lambda\kk{z}+\frac{1}{\lambda}z)}|_{\DOm},  \quad \text{for } |\lambda| = 1,
\end{align}
where the DN-map $\Lambda_{-E}$ corresponds to the potential $q_0 = 0$. In conjunction we have
\begin{align}\label{scat-trans-DN}
\vct{t}(\lambda) &= \int_{\bound}e^{\frac{i\sqrt{E}}{2}(\kk{\lambda}z+\kk{z}/\kk{\lambda})}(\Lambda_q-\Lambda_{-E})\psi(z,\lambda)|dz|,\quad \text{for } |\lambda| \neq 1 \textrm{ (non-exceptional)}, \\ \label{scat-trans-DN2}
h_{\pm}(\lambda,\lambda') &= \frac{1}{(2\pi)^2}\int_{\bound}e^{-\frac{i\sqrt{E}}{2}(\lambda'\kk{z}+z/ \lambda')}(\Lambda_q-\Lambda_{-E})\psi_{\pm}(z,\lambda)|dz|,
\end{align}
for $|\lambda| = |\lambda'|= 1$, where $|dz|$ stands for the \textit{surface measure} on $\partial \Omega$.\smallskip

Thus we have the necessary steps to reconstruct the potential $q_0$ from the DN-map $\Lambda_q$, namely:
\begin{enumerate}
\item Solve $\psi(\cdot,\lambda)|_{\DOm}$ and $\psi_{\pm}(\cdot,\lambda)|_{\DOm}$ from the boundary integral equations \eqref{BIE-lambda} and \eqref{BIE-lambdapm}, respectively.
\item Compute the scattering transforms $\vct{t}(\lambda)$ using \eqref{scat-trans-DN} and $h_{\pm}(\lambda,\lambda')$ using \eqref{scat-trans-DN2}.
\item Compute $\rho(\lambda,\lambda')$ solving one of the equations \eqref{defrho}.
\item Choose a reconstruction point $z'$ and solve $\mu(z',\lambda)$ from \eqref{IEIS-compact}.
\item Compute $q_0(z')$ from \eqref{q0-recon1}.
\end{enumerate}\smallskip

We will restrict the class of potentials to those with small (classical) fixed-energy scattering amplitude, i.e. small $\rho$. A potential such that $\rho \equiv 0$ is said to be \textit{transparent} and it is well known that there are no non-zero compactly supported transparent potentials \cite{Novikov1986b}. However, since the scattering transform $\vct{t}$ is related to $\rho$ by analytic continuation techniques (see \cite[Section 7]{Novikov1992}), its size (as well as the size of the related potential), roughly speaking, can be large even for small $\rho$. Thus, for potentials with small $\rho$ the algorithm is simplified by assuming $\rho \equiv 0$ and using only the term $\vct{t}$. We quantify the error in Lemma \ref{errrho}. See \cite{Grinevich1995} for more discussions about transparent potentials and \cite{Greenleaf2003a, Greenleaf2007, Greenleaf2009, Greenleaf2009a} for the similar phenomenon of \textit{invisibility}.
\smallskip

Thus, we will use the following algorithm:
\begin{enumerate}
\item Solve $\psi(\cdot,\lambda)|_{\DOm}$ from the boundary integral equations \eqref{BIE-lambda}.
\item Compute the scattering transform $\vct{t}(\lambda)$ using \eqref{scat-trans-DN} 
\item Choose a reconstruction point $z'$ and solve $\mu(z',\lambda)$ from \eqref{IEIS-compact} without the term $\MM$.
\item Compute an approximation $\tilde q_0(z')$ of $q_0(z')$ using \eqref{q0-recon1}.
\end{enumerate}

\begin{lemma} \label{errrho}
Let $\Omega \subset \R^2$ be a open bounded domain with $C^2$ boundary and let $q_0 \in W^{m,1}(\Omega)$, real-valued, with $\mathrm{supp}(q_0) \subset \Omega$ and $m \geq 3$. Assume that $\|q_0\|_{m,1} \leq N$ and that $E > E_1(N,\Omega)$ is sufficiently large (in particular there are no exceptional points). Let $\vct{t}$ be the scattering transform defined in \eqref{scattering-transform} and $\rho$ be the function defined in \eqref{defrho1}. Let $\tilde q_0$ be the potential obtained solving the non-local Riemann-Hilbert problem \eqref{D-bar}-\eqref{corrt} with scattering data given by $\vct{t}$ and $\rho \equiv 0$. Then there is a constant $C=C(\Omega,N,m) >0$ such that
\begin{equation}\label{estrho}
\|q_0-\tilde q_0\|_{L^{\infty}(\Omega)}\leq C(\Omega,N,m)E \|\rho\|_{L^2(T^2)},
\end{equation}
where $T = \{\lambda \in \C : |\lambda| =1\}$.
\end{lemma}

The proof of Lemma \ref{errrho} is given at the end of this section.\smallskip

\begin{remark}
We want to underline that the assumptions made on the potential and the energy are needed for a rigorous justifications of our method. Numerical results presented in Section \ref{NUMERICS} strongly suggest that our algorithm performs well at low energies for any $L^{\infty}$ potential (in absence of exceptional points). See also Subsection \ref{relatedwork} for more discussions about the energy.
\end{remark}

We will now give an interpretation of our algorithm in the Born approximation. Let $E>0$ be fixed and consider the classical scattering amplitude $f(k,l)$, $k,l \in \R^2$ with $k^2 =l^2 =E$ (see for instance \cite[(1.3a)]{Novikov1992} for a definition). From \cite[Theorem 5.2, Proposition 5.1]{Novikov1992} $f$, $h_{\pm}$ and $\rho$ are connected through integral equations and, roughly speaking, they contain the same information on a potential $q_0$. Assuming the Born approximation, i.e. $\norm{q_0}_{L^{\infty}(\Omega)} \ll E$, we have that $f(k,l) \approx \F [q_0](k-l)$, where $\F$ is the 2D Fourier transform. Thus the classical scattering amplitude $f(k,l)$ at fixed energy $E>0$, or equivalently $\rho(\lambda,\lambda')$, determines $\F [q_0](p)$ for $|p| \leq 2 \sqrt{E}$. Under the same assumptions, the non-physical scattering transform $\vct{t}(\lambda)$ determines $\F [q_0](p)$ for $|p| \geq 2 \sqrt{E}$. In particular, in Section \ref{sec:truncation} we will consider a truncated scattering transform $\vct{t}_R$, which is $0$ for $|\lambda| < 1/R$ and $|\lambda| > R$, for $R > 1$, and equal to $\vct{t}$ otherwise. In the Born approximation, $\vct{t}_R$ determines $\F [q_0](p)$ for $ 2 \sqrt{E} \leq |p| \leq \sqrt{E}\left(R+\frac 1 R \right)$. Intuitively, the proposed algorithm allows the reconstruction of the Fourier transform of a potential in the annulus $2\sqrt{E} \leq |p| \leq  \sqrt{E}\left(R+\frac 1 R \right)$: thus it provides good results for data acquired at fixed energy $E$ far from $0$ and $+\infty$.

In section \ref{sec:Dbar-numerics} we give a more detailed numerical algorithm for this reconstruction procedure.
 
\subsection{Related work}\label{relatedwork}
A reconstruction algorithm for the Gel'fand-Calder\'on problem at fixed positive energy was proposed by R. Novikov and one of the authors in \cite{Novikov01012013} (see \cite{Burov2012,Burov2013} for numerical results). This algorithm and the one presented in this paper have similar theoretical background but present several differences:
\begin{itemize}
\item In the algorithm of \cite{Novikov01012013}, only the scattering function $\rho$ (or $h_{\pm}$) is reconstructed from the Dirichlet-to-Neumann map and used in the solution of the Riemann-Hilbert problem. In the present paper we use only the scattering transform $\vct{t}$.
\item The algorithm in \cite{Novikov01012013} is Lipschitz stable with an error term depending on the energy. Our proposed algorithm is logarithmic stable, with an error quantified by Lemmas \ref{errrho} and \ref{lemmaregular}. Concerning speed, the most computationally expensive steps in the algorithm in \cite{Novikov01012013} are two 1D linear integral equations, while in our method they are a 1D linear integral equation and a $\bar \partial$-equation (2D linear integral equation). The former algorithm is then faster than the latter.
\item Both algorithms, to be rigorously justified, need the energy to be sufficiently large with respect to the $L^{\infty}$ norm of the potential. But numerical evidences show that they perform well in different energy ranges: the algorithm in \cite{Novikov01012013} at high energies, while the present one at low energies. In Figure \ref{fig:comparison} we present a numerical comparison of the two algorithms at different fixed energies: it is clear that at low energies our method provides better reconstructions. Along with the results of figures \ref{fig:comparison} and \ref{fig:Etest} we conjecture that for potentials $q_0$ such that $\|q_0\|_{L^{\infty}(\Omega)}=1$ our method performs well, when $10^{-5} \leq E  \leq 5$.
\end{itemize}

In conclusion, despite our method present several disadvantages with respect to the one proposed in \cite{Novikov01012013}, it is, as far as we know, the only known algorithm producing good reconstructions at low positive energies.\smallskip

Other reconstruction algorithms for this problem have been proposed. The fundamental paper \cite{Bukhgeim2008} provides a reconstruction method without any assumptions on the energy; this algorithm is less stable than the one presented in this paper, since it's based on properties of generalized scattering data for only large complex parameters. 

Very recently, in the preprint \cite{lakshtanov2015global} a new global reconstruction method was proposed (without any assumptions on the energy). This method is based on a generalization of the Riemann-Hilbert problem introduced in Section \ref{sec:dbar}, which is able to deal with the presence of exceptional points. 

To our knowledge, no numerical studies based on \cite{Bukhgeim2008} or \cite{lakshtanov2015global} have been presented yet.

\begin{proof}[Proof of Lemma \ref{errrho}]
Thanks to the assumptions on the potential and the energy, the non-local Riemann Hilbert problem \eqref{D-bar}-\eqref{corrt} is solved with scattering data $(\vct{t},\rho)$ and $(\vct{t},0)$ (see \cite[Theorems 6.1 and 6.2]{Novikov1992} for a proof). 

Estimate \eqref{estrho} is a direct consequence of technical results of \cite{Santacesaria2015} used to prove a stability estimate for this problem. We can repeat the arguments of \cite[Section 4]{Santacesaria2015} in order to obtain the equality \cite[identity (4.12)]{Santacesaria2015}
\begin{equation} \label{idest2}
q_0(z)-\tilde q_0(z)=2 i \sqrt{E}(A-\tilde A +B- \tilde B + C- \tilde C),
\end{equation}
where $A, \tilde A, B, \tilde B,C,\tilde C$ are constructed as $A, B, C$ in \cite[Section 4]{Santacesaria2015} with $\mathrm{sgn}(|\lambda|^2-1)\frac{\vct{t}(\lambda)}{4\pi \bar \lambda}$ instead of $r(\lambda)$, and with $\rho$ and $0$ instead of $\rho$, respectively. In \cite[Section 5]{Santacesaria2015} the right hand side of \eqref{idest2} is estimated in terms of scattering data. Since both potentials are reconstructed from the same non-physical scattering transform $\vct{t}$, we have that $A-\tilde A \equiv 0$ (this follows from the estimate after \cite[estimate (5.6)]{Santacesaria2015}). For the same reason, from the estimate after \cite[estimate (5.8)]{Santacesaria2015} we have
\begin{equation}
|B-\tilde B| \leq c(\Omega,N,m)\sqrt{E}\|\rho\|_{L^2(T^2)},
\end{equation}
and for $C - \tilde C$ the same argument gives
\begin{equation}
|C-\tilde C| \leq c(\Omega,N,m)\sqrt{E}\|\rho\|_{L^2(T^2)}.
\end{equation}
The proof follows from these estimates and identity \eqref{idest2}.
\end{proof}

\section{Computation of the Faddeev Green's function for positive energy}\label{GZETA}
As in the previous section, we will identify the plane $\R^2$ as the complex plane by writing 
$$
x = [x_1\ x_2]^T\in\R^2,\quad z = x_1+ix_2\in\C.
$$
Recall also the transformation \eqref{lambda-definition} between $\lambda$ and $\zeta$ parameters.

We need a numerical algorithm for the Faddeev Green's function $g_\lambda(z)$ for any point $\abs{z}\leq 1$ and any $\lambda$ with $\abs{\lambda}>1$, see the symmetry \eqref{tsymmetry}.

We remark that in the case $\zeta\cdot\zeta = E=0$ and $\Imm(\zeta) \neq\vct{0}$ the numerical computation of $g_\zeta(z)$ was first presented in \cite{Siltanen1999} and then used in the context of the inverse conductivity problem in \cite{Siltanen2000}. The method was later refined in \cite{Ikehata2004b}. This approach is based on implementing appropriate numerical integrations depending on the location of the evaluation point $z$. The zero-energy computation can also be simply implemented using \cite[formula (3.10)]{Boiti1987} and Matlab's built-in exponential-integral function: {\tt g = exp(-1i*z).*real(expint(-1i*z))/(2*pi)}.

The Faddeev Green's function $g_\lambda(z)$ for positive energy has similar scaling and rotational properties as in the zero energy case. However, our approach is based on $g_\zeta(z)$, $\zeta\cdot\zeta=E>0$: the main reason is that we can use residue calculus to compute simpler formulas in a similar way to \cite{Siltanen1999}. When using the $\zeta$ parameters, the scaling and rotational properties will be employed differently from the zero-energy case. Thus the formulas obtained in this section cannot be directly obtained from \cite{Siltanen1999,Siltanen2000,Ikehata2004b} and this region-based method has to be modified. 

Recall the formula \eqref{green_zeta},
\begin{equation*}
g_\zeta(z) = \frac{1}{4\pi^2}\int_{\R^2}\frac{e^{iy\cdot z}}{y\cdot y+2\zeta\cdot y}dy.
\end{equation*}
The following relations can be seen.
\begin{lemma}
Let $\alpha\in \R\setminus\{0\}$, $R$ a rotational matrix with $\det(R) = 1$ and $R\zeta = R(\Ree(\zeta))+iR(\Imm(\zeta))$. Then the Faddeev Green's function $g_\zeta(z)$ with $\zeta\cdot\zeta=E>0$ satisfies
\begin{eqnarray}
g_\zeta(\alpha z) &=& g_{\alpha\zeta}(z)\label{alpha_relation} \\
g_\zeta(Rz) &=& g_{R^{-1}\zeta}(z)\label{rotation_relation} \\
g_\zeta(\matr{-x_1\\x_2}) &=& g_{\matr{-\zeta_1\\ \zeta_2}}(z)\label{x1_relation} \\
g_\zeta(\matr{x_1\\-x_2}) &=& g_{\matr{\zeta_1\\ -\zeta_2}}(z)\label{x2_relation} \\
g_{\kk{\zeta}}(z) &=& \kk{g_{\zeta}(-z)}.\label{conjugate_relation}
\end{eqnarray}
\end{lemma}
We can use the rotation relation \eqref{rotation_relation} to reduce $\zeta$ to the form
\begin{equation}\label{reduced-zeta}
\zeta = \matr{k_1\\0}+i\matr{0\\k_2},\quad \abs{k_1}>k_2>0.
\end{equation}
For this reduced $\zeta$, using relations \eqref{conjugate_relation} and \eqref{x2_relation} we have
\begin{equation*}
\kk{g_\zeta(-z)} = g_{\kk{\zeta}}(z) = g_\zeta(\matr{x_1\\-x_2}),
\end{equation*}
which results to a switching relation
\begin{equation}\label{x1-symmetry}
g_\zeta(\matr{-x_1\\x_2}) = \kk{g_\zeta(\matr{x_1\\x_2})}
\end{equation}
The strategy for computing $g_\lambda(z)$ is now the following:
\begin{enumerate}
\item Use \eqref{lambda-definition} to compute $\zeta$ from $\lambda$.
\item Find the rotational matrix $R$ that satisfies $R (\Imm(\zeta)) = [0,k_2]^{T}$ for some $k_2>0$; then write
  \begin{equation*}
    g_{\zeta}(z) = g_{R^{-1}R\zeta}(z) = g_{R\zeta}(Rz),
  \end{equation*}
  where $R\zeta$ is in the reduced form \eqref{reduced-zeta}.
\item The smaller $\abs{z}$ the more computational problems we have as will be seen later. Thus, for very small $\abs{z}$ we use a method of single layer potential described in section \ref{singlelayer}. Further, use relation \eqref{alpha_relation} to scale points outwards and relation \eqref{x1-symmetry} to switch from $x_1<0$ to $x_1\geq0$.
\item Use computational domains to compute $g_{\zeta}(z)$ for reduced $\zeta$ and $z$ with $x_1\geq 0$.
\end{enumerate}
It takes some analysis to find suitable computational domains for the last step.

Assume we have the reduced $\zeta$ of \eqref{reduced-zeta} and the switched $z=x_1+ix_2$ with $x_1\geq0$. Write $t = y_1+k_1$, $a = (y_2+k_2i)^2-E$ and subsequently the denominator of the integrand in \eqref{green_zeta} as
\begin{eqnarray*}
y\cdot y + 2y\cdot\zeta &=& y_1^2 + y_2^2 + 2y_1k_1+2y_2k_2i \\
&=& (y_1+k_1)^2+(y_2+k_2i)^2-E \\
&=& t^2+a \\
&=& (t+i\sqrt{a})(t-i\sqrt{a})
\end{eqnarray*}
We define the square root in the same way MATLAB calculates it by default, that is for a complex number $z=r\exp(i\theta),0\leq\theta<2\pi,r\geq0$ the square root is
\begin{equation*}
\sqrt{z} = \left\{
\begin{array}{cl}
\sqrt{r}e^{i\theta/2}&,\quad 0\leq\theta\leq\pi,\\
\sqrt{r}e^{-i(2\pi-\theta)/2}&,\quad \pi<\theta<2\pi.
\end{array}
\right.
\end{equation*}
This way the square root has the following properties: for any $z\in\C$ we have
\begin{eqnarray}
  \Ree(\sqrt{z}) &\geq& 0 \label{sqrt-Re}\\
\sqrt{\kk{z}} &=& \kk{\sqrt{z}} \label{sqrt-kk}.
\end{eqnarray}
The numerator of the integrand in \eqref{green_zeta} becomes
\begin{eqnarray*}
e^{iz\cdot y} &=& e^{i(x_1(t-k_1)+x_2y_2)} \\
&=& e^{i(x_2y_2-x_1k_1)}e^{ix_1t}.
\end{eqnarray*}
The integral in \eqref{green_zeta} is thus transformed into
\begin{equation*}
\int_{\R^2}\frac{e^{iy\cdot z}}{y\cdot y+2\zeta\cdot y}dy = \int_{-\infty}^{\infty}e^{i(x_2y_2-x_1k_1)}\left(\int_{-\infty}^{\infty}\frac{e^{ix_1t}}{(t+i\sqrt{a})(t-i\sqrt{a})}dt\right)dy_2.
\end{equation*}
The integral over the real parameter $t$ is complexified with $w=w_R+iw_I\in\C$ and we write
\begin{eqnarray}\label{fw}
f(w) = \frac{e^{ix_1w}}{(w+i\sqrt{a})(w-i\sqrt{a})} = \frac{e^{ix_1w_R}e^{-x_1w_I}}{(w+i\sqrt{a})(w-i\sqrt{a})},
\end{eqnarray}
The poles of the function $f(w)$ are $\pm i\sqrt{a}$, where 
\begin{equation*}
a = (y_2+k_2i)^2-E = y_2^2-E-k_2^2 + 2y_2k_2i. 
\end{equation*}
It follows from our definition of $\sqrt{\cdot}$ that
\begin{itemize}
\item when $y_2>0$, $a$ is in the upper half plane, so $i\sqrt{a}$ is in quadrant 2 and $-i\sqrt{a}$ in quadrant 4 (note, that $k_2>0$), and
\item  when $y_2<0$, $a$ is in the lower half plane, so $i\sqrt{a}$ is in quadrant 1 and $-i\sqrt{a}$ in quadrant 3.
\end{itemize}
When $w_I,x_1\geq0$ we have
\begin{equation*}
\abs{f(w)}\rightarrow 0,\quad \mbox{ as }\abs{w}\rightarrow\infty.
\end{equation*}
We choose the integration path 
\begin{equation*}
\Gamma = [-R,R]\cup\{R\exp(i\theta): 0\leq\theta\leq\pi\},
\end{equation*}
so when  $R\rightarrow\infty$ the pole $w=i\sqrt{a}$ is inside the path. Using residue calculus we get
\begin{eqnarray*}
\int_{\R\subset\C}f(w)dw &=& \int_{\Gamma}f(w)dw = 2\pi i\res_{w=i\sqrt{a}}f(w) \\
&=& 2\pi i\lim_{w\rightarrow i\sqrt{a}}(w-i\sqrt{a})f(w) = 2\pi i\frac{e^{ix_1(i\sqrt{a})}}{i\sqrt{a}+i\sqrt{a}}\\
&=& \pi\frac{e^{-x_1\sqrt{a}}}{\sqrt{a}}
\end{eqnarray*}
Thus
\begin{eqnarray*}
g_\zeta(z) &=& \frac{1}{(2\pi)^2}(\int_{-\infty}^{0}e^{i(x_2y_2-x_1k_1)}\pi \frac{e^{-x_1\sqrt{a}}}{\sqrt{a}} dy_2 + \int_{0}^{\infty}e^{i(x_2y_2-x_1k_1)}\pi \frac{e^{-x_1\sqrt{a}}}{\sqrt{a}}dy_2) \\
&=& \frac{1}{4\pi}e^{-ix_1k_1}(\int_{0}^{\infty}e^{-ix_2y_2}\frac{e^{-x_1\sqrt{\kk{a}}}}{\sqrt{\kk{a}}} dy_2 + \int_{0}^{\infty}e^{ix_2y_2}\frac{e^{-x_1\sqrt{a}}}{\sqrt{a}}dy_2).
\end{eqnarray*}
Because of \eqref{sqrt-Re} we have
\begin{equation*}
e^{-ix_2y_2}\frac{e^{-x_1\sqrt{\kk{a}}}}{\sqrt{\kk{a}}} = \kk{e^{ix_2y_2}\frac{e^{-x_1\sqrt{a}}}{\sqrt{a}}},
\end{equation*}
and so the following formula is obtained.
\begin{lemma}
For $x_1\geq0$,
\begin{equation}\label{gz1}
g_\zeta(z) = \frac{1}{2\pi}e^{-ix_1k_1}\Ree(\int_{0}^{\infty}e^{ix_2t}\frac{e^{-x_1\sqrt{(t+k_2i)^2-E}}}{\sqrt{(t+k_2i)^2-E}}dt).
\end{equation}
\end{lemma}
We want to numerically compute the integral of \eqref{gz1}. We can only compute up to a finite limit, say from $0$ to $T$. Two problems might occur, the integrand either converges slowly, meaning we have to take $T$ very large, or the integrand might oscillate fast, meaning we have to take a great number of integration points in $[0,T]$. We see that problematic situations in using \eqref{gz1} arise when $\abs{x_2}$ is large (oscillation), $x_1$ is small (convergence) and when $k_2$ is large (oscillation). When $x_1$ is large we have oscillation, but also better convergence, meaning \eqref{gz1} is usable.

Also worth noting is that when $\lambda$ is close to one, then $k_2$ is close to zero and the poles $\pm i\sqrt{a}$ are close to the integration path $\Gamma$ causing numerical problems.

These observations lead to additional versions of formula \eqref{gz1}, used by the different computational domains. Write
\begin{equation*}
  g(w) = e^{ix_2 w}\frac{e^{-x_1\sqrt{a}}}{\sqrt{a}},\quad a = (w+k_2i)^2-E,\quad w = w_1+iw_2\in\C,
\end{equation*}
and consider the complexified integral $\int_{\R^+}g(w)dw$ of \eqref{gz1}. We have
\begin{equation*}
  e^{ix_2 w} = e^{ix_2w_1}e^{-x_2w_2},
\end{equation*}
so when $x_2,w_2\geq0$ or $x_2,w_2<0$ we have 
\begin{equation}\label{gz-convergence}
\abs{g(w)}\rightarrow0\mbox{ as }\abs{w}\rightarrow\infty. 
\end{equation}
This is because the numerator $\exp(-x_1\sqrt{a})$ converges to zero, since $x_1\geq0$ and \eqref{sqrt-Re} holds. The branch points of $g(w)$ are $\pm\sqrt{E}-k_2i$. 
\begin{lemma}
Let $x_1\geq 0$ and $x_2\geq 0$, then
\begin{eqnarray}
g_\zeta(z) &=& \frac{1}{2\pi}e^{-ix_1k_1}\Ree(\int_0^\infty e^{-x_2 t}\frac{e^{-x_1 i\sqrt{t^2+k_1^2+2tk_2}}}{\sqrt{t^2+k_1^2+2tk_2}}dt)\label{gz2}.
\end{eqnarray}
\begin{proof}
We choose the integration path
\begin{equation*}
\Gamma_1 = [0,R]\cup\{R\exp(i\theta): 0\leq\theta\leq\pi/2\}\cup\{iR(1-t):  0\leq t\leq 1\}.
\end{equation*}
We have $\int_{\Gamma_1}g(w)dw = 0$, since $g(w)$ is analytic in the first quadrant. Because of \eqref{gz-convergence} we then have
\begin{equation*}
  \int_{\R^+}g(w)dw = -\int_{\infty}^0g(it)idt = \int_0^\infty e^{-x_2 t}\frac{e^{-x_1 i\sqrt{(t+k_2)^2+E}}}{\sqrt{(t+k_2)^2+E}}dt. \qedhere
\end{equation*}
\end{proof}
\end{lemma}
The integrand in \eqref{gz2} converges to zero quickly for large $x_2$ and has high oscillation for large $x_1$.
\begin{lemma}
Let $x_1\geq 0$ and $x_2<0$, then
\begin{eqnarray}
g_\zeta(z) &=& \frac{1}{2\pi}e^{-ix_1k_1}\Ree(I_1-ie^{i(\sqrt{E}+1)x_2}\int_0^\infty e^{x_2 t}\frac{e^{-x_1\sqrt{b}}}{\sqrt{b}}dt)\label{gz3},
\end{eqnarray}
where $I_1 = \int_0^{\sqrt{E}+1}g(t)dt$ and $b=(\sqrt{E}+1+(k_2-t)i)^2-E$.
\begin{proof}
Define the paths
\begin{eqnarray*}
P_1 &=& [0,\sqrt{E}+1],\\
L_1 &=& \{\sqrt{E}+1-iRt: 0\leq t\leq 1\},\\
L_2 &=& \{\sqrt{E}+1+R\exp(i\theta): \frac{3}{2}\pi\leq\theta\leq 2\pi\},\\
L_3 &=& \{\sqrt{E}+1+(1-t)R: 0\leq t\leq 1\}.
\end{eqnarray*}
The branch point $\sqrt{E}-k_2i$ is avoided by integrating along $\Gamma_2 = P_1\cup (L_1\cup L_2\cup L_3)$, where $g(w)$ is analytical inside the loop $L_1\cup L_2\cup L_3$, and $\abs{g(w)}\rightarrow 0$ on the circle $L_2$ as $R\rightarrow\infty$. Thus
\begin{equation*}
  \int_{\R^+}g(w)dw = I_1 + \int_0^\infty g(\sqrt{E}+1-it)(-i)dt,
\end{equation*}
where
\begin{equation*}
 \int_0^\infty g(\sqrt{E}+1-it)(-i)dt = \int_0^\infty e^{i(\sqrt{E}+1)x_2}e^{x_2 t}\frac{e^{-x_1\sqrt{b}}}{\sqrt{b}}(-i)dt. \qedhere
\end{equation*}
\end{proof}
\end{lemma}
The integrand in \eqref{gz3} converges to zero quickly for large $\abs{x_2}$ and has high oscillation for large $x_1$.

\subsection{Choosing upper limits for the integrals}
Write $g_\zeta^{T_1}$,$g_\zeta^{T_2}$ and $g_\zeta^{T_3}$ for the finite integrals for \eqref{gz1},\eqref{gz2} and \eqref{gz3} respectively. We need to choose the upper limits $T_i$, $i=1,2,3$. There will be numerical error caused by the neglected part of the integral and the numerical integration method used.  It is decided to simply require
\begin{equation}\label{error-requirement}
  \abs{g_\zeta-g_\zeta^{T_i}} < 10^{-8}.
\end{equation}
The error induced by the numerical integration method is assumed not to be dependant on $\lambda$ or $z$. For $g_\zeta^{T_i}$, the integration range $[0,T_i]$ is divided into $M_i$ points (with $g_\zeta^{T_3}$ there is also the additional integral $I_1$) and the Gaussian quadrature is used. The integers $M_i$ are chosen large enough so that for any integer $M>M_i$ the first 8 digits are not changing in the numerical value of $g_\zeta^{T_i}(z)$. In this test the choice of $z$ has only a minor effect, it is done by choosing the ``worst'' possible point for any given computational domain; for example, for $g_\zeta^{T_1}(z)$ using $z = [1,1]$ the integrand has more oscillation than with the point $z=[1,0]$, and thus needs a larger parameter $M_1$.

Finding $T_i$ is a bit cumbersome. The following proposition guarantees us the error requirement \eqref{error-requirement}.
\begin{Prop}
Choose
\begin{eqnarray}
T_1 &=& \max\{\frac{14\cdot 2^{1/4}}{x_1c_1},2k_1\}\label{T1-def},\\
T_2 &=& \frac{14}{x_2}\label{T2-def},\\
T_3 &=& \frac{14}{c_2x_1-x_2}+k_2\label{T3-def},
\end{eqnarray}
where $c_1,c_2$ are constants depending on $k_1$ and $k_2$. Then,
\begin{equation*}
  \abs{g_\zeta(z)-g_\zeta^{T_i}(z)} <10^{-8},\quad i=1,2,3.
\end{equation*}
\end{Prop}

\begin{proof}
In \eqref{gz1} we have the term $a$ and
\begin{eqnarray*}
  \abs{\sqrt{a}} &=& \abs{\sqrt{(t^2-k_2^2-E)+2k_2ti}} = ((t^2-k_1^2)^2+4k_2^2t^2)^{1/4}\nonumber\\
&=& (t^4-2t^2k_1^2+k_1^4+4k_2^2t^2)^{1/4}\nonumber\\
&\geq& (t^4-2k_1^2t^2)^{1/4} \geq (t^4-1/2t^4)^{1/4}\nonumber\\
&=& \frac{t}{2^{1/4}} \label{gz1-ul},
\end{eqnarray*}
when $t\geq 2\abs{k_1}$. Then, writing $\theta$ for the angle $\sqrt{a} = r\exp(i\theta)$,
\begin{eqnarray*}
  \abs{\frac{e^{-x_1\sqrt{a}}}{\sqrt{a}}} &\leq& \frac{e^{-x_1\Re(\sqrt{a})}}{t/2^{1/4}} = \frac{e^{-x_1\cos(\theta)\abs{\sqrt{a}}}}{t/2^{1/4}} \\
&\leq& \frac{e^{-x_1\cos(\theta)t/2^{1/4}}}{t/2^{1/4}}.
\end{eqnarray*}
The angle goes to zero as $t\rightarrow\infty$, so $\cos(\theta)\rightarrow 1^-$. Since $t\geq 2\abs{k_1}$ we write $c_1=\cos(\theta_1)$, where the angle of $\sqrt{a}|_{t=2k_1}$ is $\theta_1$, and so we have for the integral
\begin{eqnarray}
 \abs{\int_{T}^\infty e^{ix_2t}\frac{e^{-x_1\sqrt{a}}}{\sqrt{a}}dt} &\leq& \int_{T}^\infty \frac{e^{-x_1c_1t/2^{1/4}}}{t/2^{1/4}}dt \nonumber\\
&=& 2^{1/4}\int_{x_1c_1 T/2^{1/4}}^\infty\frac{e^{-s}}{s}ds = 2^{1/4}E_i(x_1c_1 T/2^{1/4}) \label{remainder}.
\end{eqnarray}
The exponential integral function $E_i$ can be computed in MATLAB with {\tt expint.m}. Because of \eqref{error-requirement} we require that the remainder \eqref{remainder} is of the order $2\pi/2^{1/4}\cdot 10^{-8} \approx 7.47\cdot 10^{-8}$. We can test with MATLAB that $E_i(14)<6\cdot 10^{-8}$, so we get
\begin{eqnarray*}
\begin{array}{rcl}
  x_1c_1T/2^{1/4} &=& 14
\end{array}
\end{eqnarray*}
from which \eqref{T1-def} follows.

From \eqref{gz2} we easily get
\begin{eqnarray*}
  \abs{\int_{T_2}^\infty e^{-x_2t}\frac{e^{-x_1i\sqrt{t^2+k_1^2+2tk_2}}}{\sqrt{t^2+k_1^2+2tk_2}}dt} &\leq& \int_{T_2}^\infty\frac{e^{-x_2 t}}{t}dt \\
&=& E_i(x_2T_2).
\end{eqnarray*}
Thus the upper limit \eqref{T2-def} follows, as before, from
\begin{eqnarray*}
  x_2T_2 &=& 14.
\end{eqnarray*}

Starting from \eqref{gz3} we have
\begin{eqnarray*}
  \abs{\sqrt{b}} &=& ((k_2-t)^4-2(2\sqrt{E}+1)(k_2-t)^2+(2\sqrt{E}+1)^2\\
& &+4(\sqrt{E}+1)^2(k_2-t)^2)^{1/4}\\
&\geq& ((k_2-t)^4)^{1/4} = \abs{t-k_2}.
\end{eqnarray*}
Using the same argument as preceding \eqref{remainder}, we write $c_2 = \cos(\theta_2)$, where the angle of $\sqrt{b}|_{t=k_2}$ is $\theta_2$. Then for $T_3> k_2$ we have
\begin{eqnarray*}
 \abs{\int_{T_3}^\infty e^{x_2t}\frac{e^{-x_1\sqrt{b}}}{\sqrt{b}}dt} &\leq& \int_{T_3}^\infty \frac{e^{x_2t-x_1c_2(t-k_2)}}{t-k_2}dt \\
&=& e^{x_2k_2}\int_{(c_2x_1-x_2)(T_3-k_2)}^\infty\frac{e^{-s}}{s}ds \\
&=& e^{x_2k_2}E_i((c_2x_1-x_2)(T_3-k_2))\\
&\leq& E_i((c_2x_1-x_2)(T_3-k_2)),
\end{eqnarray*}
and \eqref{T3-def} follows from
\begin{equation*}
  (c_2x_1-x_2)(T_3-k_2) = 14. \qedhere
\end{equation*}
\end{proof}

\subsection{Use of single-layer potential for small $z$}\label{singlelayer}
For small values of $z$ there is a problem of slow convergence. We will evade this problem by the use of the single-layer potential for a function that satisfies the Helmholtz equation. Write 
$$
E = k^2,\quad G_\zeta(z) = \exp(i\zeta\cdot z)g_\zeta(z),\quad G(z) = iH_0^{1}(k\abs{z})/4,
$$
where $H_0^1$ is Hankel's function of the first type. We have
\begin{eqnarray*}
(-\triangle-k^2)G_\zeta(z) &=& \delta_z \\
(-\triangle-k^2)G(z) &=& \delta_z ,
\end{eqnarray*}
so
\begin{equation*}
(-\triangle-k^2)(G_\zeta - G) = 0.
\end{equation*}
Write $H_\zeta:=G_\zeta-G$. For any radius $R$ there exists a single-layer potential $p(z)$, which gives the value of $H_\zeta$ by the integral
\begin{equation}\label{layer-potential}
H_\zeta(z) = \int_{\partial D(0,R)}\frac{i}{4}H_0^{1}(k\abs{z-y})p(y)d\mu(y) := S(p(\cdot))(z).
\end{equation}
Assume we know $H_\zeta(z)$ on the circle $\partial D(0,R)$, where $R$ is large enough so that we don't have the problems of slow convergence. The potential can be recovered by the inverse of the integral operator, $p = S^{-1}(H_\zeta(z))$. Then $H_\zeta(z)$ can be calculated using \eqref{layer-potential} for any $\abs{z}<R$. Finally we have
\begin{equation}\label{gzeta-layerpotential}
g_\zeta(z) = e^{-i\zeta\cdot z}(H_\zeta(z) + G(z)).
\end{equation}
The numerical implementation of this submethod is straightforward with the additional trick that the potential $p(z)$ is computed on a circle of radius $R+\epsilon>R$ so that we avoid singularities in the operator $S$.


\subsection{Computational domains and the computation of $g_\zeta(z)$}

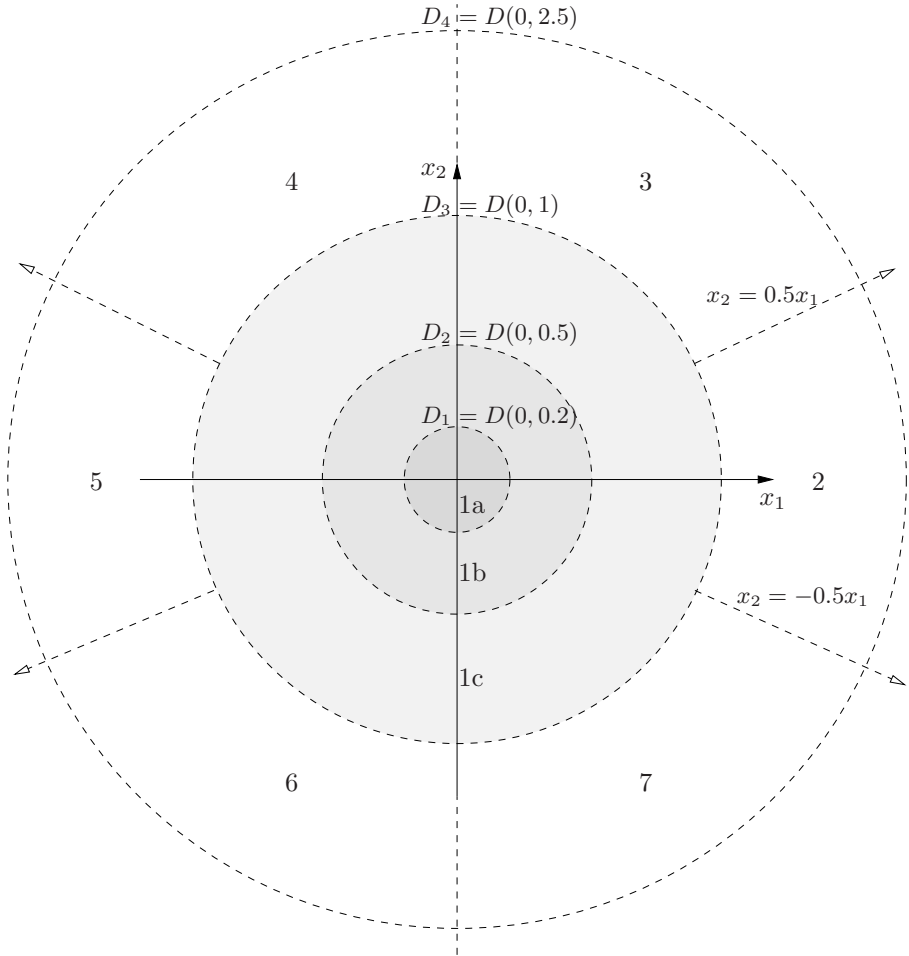
\begin{figure}[h]
\unitlength=1mm
\begin{picture}(120,140)
\epsfxsize=12cm
\put(5,5){\epsffile{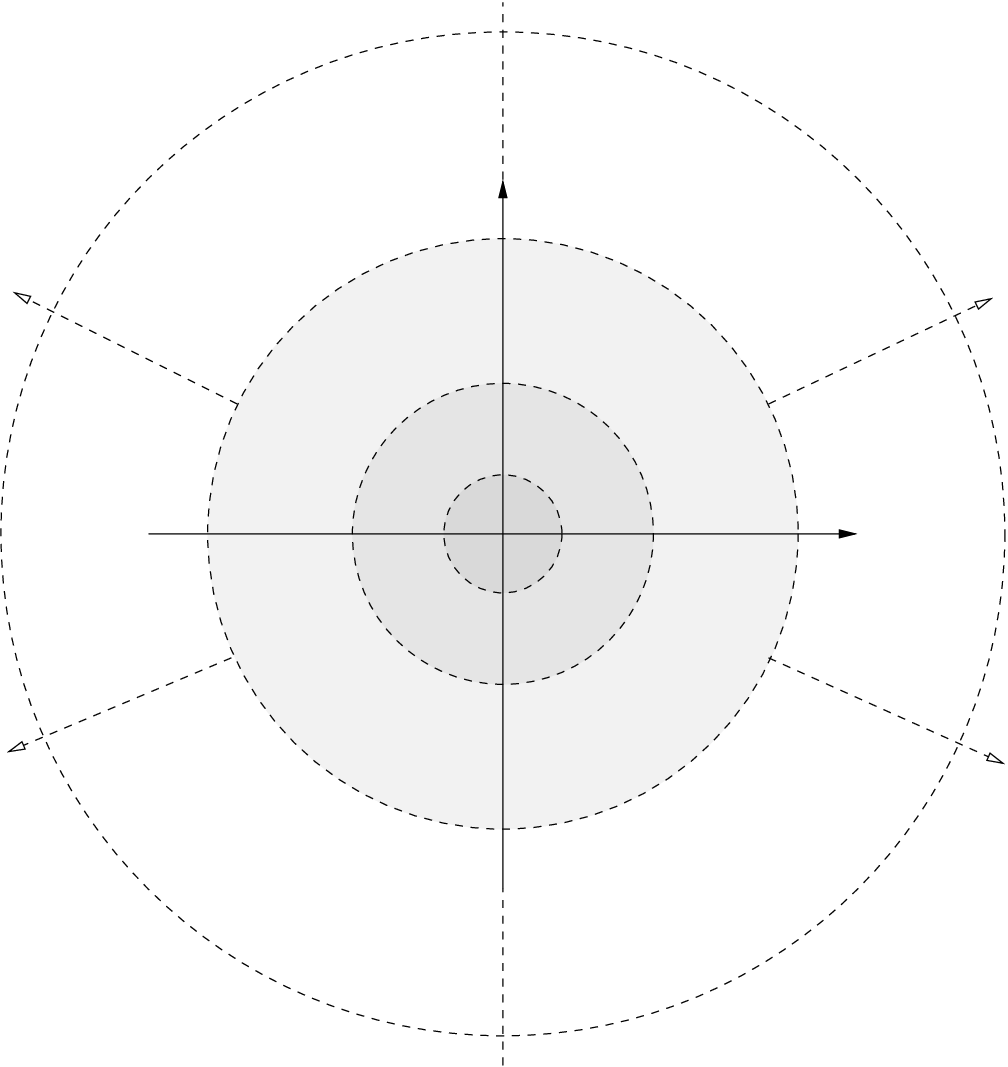}}
\put(60,109){$x_2$}
\put(105,65){$x_1$}
\put(65,64){1a}
\put(65,55){1b}
\put(65,41){1c}
\put(112,67){2}
\put(16,67){5}
\put(89,107){3}
\put(89,27){7}
\put(42,107){4}
\put(42,27){6}
\put(60,76){\small $D_1 = D(0,0.2)$}
\put(60,87){\small $D_2 = D(0,0.5)$}
\put(60,104){\small $D_3 = D(0,1)$}
\put(60,129){\small $D_4 = D(0,2.5)$}
\put(98,92){\small $x_2 = 0.5 x_1$}
\put(102,52){\small $x_2 = -0.5 x_1$}
\end{picture}
\caption{\label{areas}Computational domains. Domain 1a is the disk $D_1$, domain 1b is the annulus $D_2\setminus D_1$, domain 1c is the annulus $D_3\setminus D_2$. Domains 2,3,4,5,6 and 7 form the annulus $D_4\setminus D_3$.}
\end{figure}

\begin{figure}
\unitlength=1mm
\begin{picture}(120,60)
\epsfxsize=6cm
\put(0,0){\epsffile{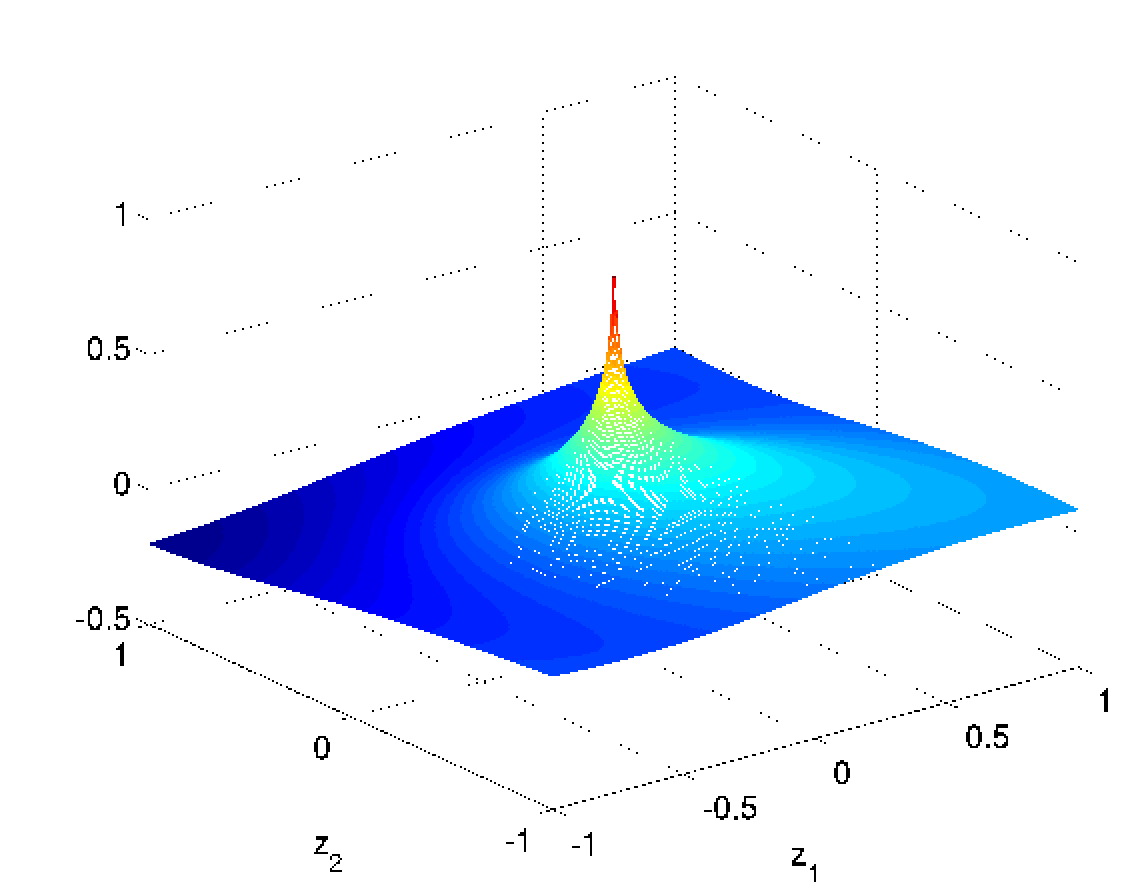}}
\epsfxsize=6cm
\put(62,0){\epsffile{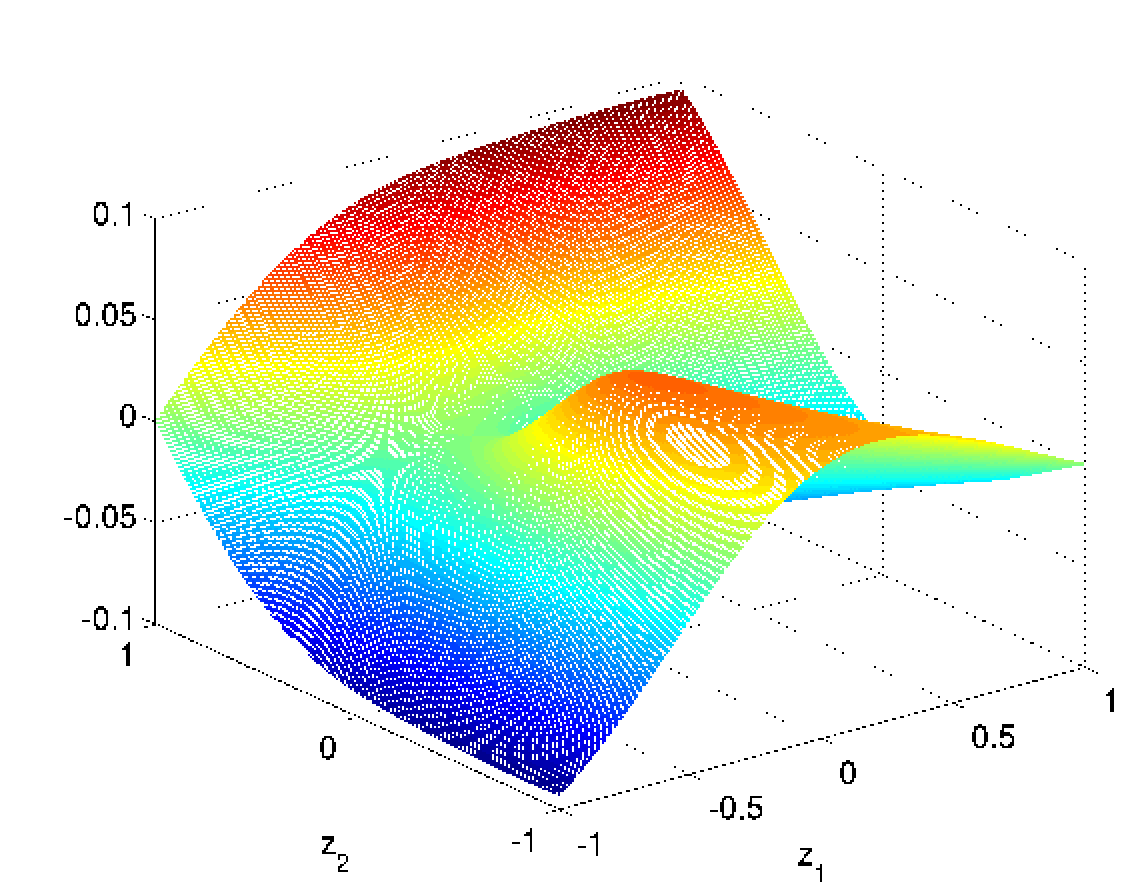}}
\put(90,49){\small $\Imm(g_\lambda(z))$}
\put(28,49){\small $\Ree(g_\lambda(z))$}
\end{picture}
\caption{\label{fig:gzeta}The real and imaginary parts of $g_{\lambda}(z)$ in $400\times400$ grid of points $z$, $\lambda=1+i$, $E=1$.}
\end{figure}

For the reduced $\zeta$ we now have the equations
\begin{equation}
g_\zeta^{T_1}(z) = \frac{1}{2\pi}e^{-ix_1k_1}\Ree(\int_{0}^{T_1}e^{ix_2t}\frac{e^{-x_1\sqrt{t^2+2tk_2i-k_1^2}}}{\sqrt{t^2+2tk_2i-k_1^2}}dt),\label{gz-T1}
\end{equation}
\begin{equation}
g_\zeta^{T_2}(z) = \frac{1}{2\pi}e^{-ix_1k_1}\Ree(\int_0^{T_2} e^{-x_2 t}\frac{e^{-x_1 i\sqrt{t^2+2tk_2+k_1^2}}}{\sqrt{t^2+2tk_2+k_1^2}}dt)\label{gz-T2},\quad x_2\geq0,
\end{equation}
\begin{eqnarray}
g_\zeta^{T_3}(z) &=& \frac{1}{2\pi}e^{-ix_1k_1}\Ree(\int_0^{\sqrt{E}+1}e^{ix_2t}\frac{e^{-x_1\sqrt{t^2+2tk_2i-k_1^2}}}{\sqrt{t^2+2tk_2i-k_1^2}}dt\label{gz-T3}\\
& &-ie^{i(\sqrt{E}+1)x_2}\int_0^{T_3}e^{x_2 t}\frac{e^{-x_1\sqrt{b}}}{\sqrt{b}}dt),\quad x_2<0,\nonumber
\end{eqnarray}
where $b=(\sqrt{E}+1+(k_2-t)i)^2-E$. See figure \ref{areas}. In general the point $z$ lies in one of the computational domains:
\begin{itemize}
\item In domain 1a, we use the single-layer potential and the equation \eqref{gzeta-layerpotential}.
\item In domain 1b, we scale the point $z$ to the annulus $D_4\setminus D_3$ using \eqref{alpha_relation}, 
  \begin{equation*}
    g_\zeta(z) = g_\zeta(1/5\cdot5x) = g_{\zeta/5}(5x)
  \end{equation*}
(note that the energy $E$ changes via this scaling transformation).
\item In domain 1c we do the same as above with the scaling factor 2.
\item In domain 2 we use \eqref{gz-T1}, since $\abs{x_2}$ is small and $x_1$ is large. The upper limit $T_1$ is computed from \eqref{T1-def}.
\item In domain 3 we use \eqref{gz-T2}, since $x_1$ is small and $x_2>0$ is large. The upper limit $T_2$ is computed from \eqref{T2-def}.
\item In domain 7 we use \eqref{gz-T3}, since $x_1$ is small, $x_2<0$ and $\abs{x_2}$ is large. The upper limit $T_3$ is computed from \eqref{T3-def}.
\item In domains 4,5,6 we use \eqref{x1-symmetry} to switch them to domains 3,2,7 respectively.
\end{itemize}
A sample of the function $g_\lambda(z)$ is pictured in \ref{fig:gzeta}, in $400\times400$ -grid of points $z$, $\lambda=1+i$, $E=1$.

\section{Numerical implementation of the D-bar method}\label{sec:Dbar-numerics}
\subsection{Simulation of measurement data}
We use truncated Fourier basis to approximate operators on the boundary $\DOm$ of the unit disk $\Omega = D(0,1)$ by finite matrices. Choose an integer $N>0$ and define the following basis functions:
\begin{equation}\label{basisfunctions}
\phi^{(n)}(\theta) = \frac{1}{\sqrt{2\pi}}e^{in\theta},\quad
n=-N,...,N.
\end{equation}
The data of the inverse problem is the DN-map \eqref{DNq-def}. Solve the problem
\begin{equation}\label{gelfand}
  (-\Delta+q) u^{(n)} = 0\mbox{ in }\Omega,\qquad
  u^{(n)} = \phi^{(n)} \mbox{ on }\DOm
\end{equation}
for $u^{(n)}$ using Finite Element Method. Define the matrix $\mathrm{L}_q=[\widehat{u}(\ell, n)]$ by
\begin{equation}\label{Lq-definition}
   \widehat{u}(\ell, n) = \int_{\DOm} \frac{\partial u^{(n)}}{\partial\nu} \overline{\phi^{(\ell)}} ds.
\end{equation}
Here $\ell$ is the row index and $n$ is the column index. The integration can be computed when the set $[0,2\pi)$ is divided into discrete points. The matrix $\mathrm{L}_q$ represents the operator $\Lambda_q$ approximately.

We add simulated measurement noise by defining
\begin{equation}\label{noiselevel}
  \mathrm{L}_q^\epsilon := \mathrm{L}_q + c\cdot\mathrm{G},
\end{equation}
where $\mathrm{G}$ is a $(2N+1)\times (2N+1)$ matrix with random entries
independently distributed according to the Gaussian normal density
$\mathcal{N}(0,1)$. The constant $c>0$ can be adjusted for different relative errors $\|\mathrm{L}_q^\epsilon-\mathrm{L}_q\|/\norm{\mathrm{L}_q}$, where $\norm{\cdot}$ is the standard matrix norm.

The DN-map $\Lambda_{-E}$ is represented by the matrix $\mathrm{L}_{-E}$ in a similar way, in the boundary value problem \eqref{gelfand} we then have $q=-E$.

\subsection{Solving the boundary integral equation}
For any $\lambda$ with $\abs{\lambda}\neq 1$ we can compute the matrix representation $\mathrm{S}_\lambda=[\widehat{s}(\ell, n)]$ of the operator $\SSS_\lambda$ \eqref{Slambda-def} defined by
\begin{equation}\label{Slambda-matrix}
   \widehat{s}(\ell, n) = \int_{\DOm} s^{(n)} \overline{\phi^{(\ell)}} ds,\quad s^{(n)} = \int_{\DOm}G_\lambda(z-y)\phi^{(n)}(y)ds(y).
\end{equation}
Here $\ell$ is the row index and $n$ is the column index, and the set $[0,2\pi)$ is divided into discrete points. The functions 
$$
\psi(z,\lambda)|_{\DOm},\quad \exp(i\sqrt{E}/2(\lambda\kk{z}+z/\lambda))|_{\DOm},
$$
can be expressed as vectors $\psi_\lambda^{\textrm{vec}}, e_\lambda^{\textrm{vec}}$ respectively using the basis functions $\phi^{(n)}$. Then, the boundary integral equation \eqref{BIE-lambda} is approximated by the equation
\begin{align}\label{BIE-discrete}
(\textrm{I}+\mathrm{S}_\lambda(\mathrm{L}_q-\mathrm{L}_{-E}))\psi_\lambda^{\textrm{vec}} = e_\lambda^{\textrm{vec}},
\end{align}
where $\textrm{I}$ is the correct sized unit matrix. These are easily solved for the vectors $\psi_\lambda^{\textrm{vec}}$, by inverting the matrix $\textrm{I}+\mathrm{S}_\lambda(\mathrm{L}_q-\mathrm{L}_{-E})$.

\subsection{Truncation of the scattering transform}\label{sec:truncation}
The computation of the CGO solutions for $|\lambda|$ close to $0$ and $+\infty$ is computationally unstable. For this reason we will calculate the values of the scattering transform only when $1/R < |\lambda| < R$, for some $R > 1$ fixed. The following lemma rigorously justifies the use of such a truncation and gives an explicit estimate, assuming some smoothness of the potential. As a corollary we obtain that the low frequency part of the potential is asymptotically close to its \textit{non-linear low frequency} part. More precisely, the potential reconstructed from the truncated scattering transform on the annulus $\{1/R < |\lambda| < R \}$ and the one obtained from the truncated Fourier transform on a ball of radius $\sqrt{E}R$, coincide up to $O((\sqrt{E}R)^{-(m-2)})$ where $m$ is related to the regularity of the potential. See Corollary \ref{corolregular} for clarity.

\begin{lemma}\label{lemmaregular}
Let $\Omega \subset \R^2$ be a open bounded domain with $C^2$ boundary and let $q_0 \in W^{m,1}(\Omega)$, real-valued, with $\mathrm{supp}(q_0) \subset \Omega$ and $m \geq 3$. Assume that $\|q_0\|_{m,1} \leq N$ and that $E > E_1(N,\Omega)$ is sufficiently large, so that there are no exceptional points. Fix $R \geq R_0(N,m) >2$ and define $\vct{t}_R (\lambda) = \vct{t}(\lambda)\chi_R(\lambda)$, where $\vct{t}(\lambda)$ is defined in \eqref{scattering-transform} and $\chi_R$ is the characteristic function of the annulus $\A_R =\{ \lambda \in \C : 1/R \leq \lambda \leq R\}$. Let $\rho$ be the function defined in \eqref{defrho1}. Let $q_R$ be the potential obtained solving the non-local Riemann-Hilbert problem \eqref{D-bar}-\eqref{corrt} with scattering data given by $\vct{t}_R$ and $\rho$. Then there is a constant $C=C(\Omega,N,m) >0$ such that
\begin{equation}\label{regular}
\|q_0 -q_R\|_{L^{\infty}(\Omega)} \leq C E^{-(m-2)/2}R^{-(m-2)}.
\end{equation}
\end{lemma}

\begin{proof}
First we must verify that the potential $q_R$ is well defined, that is that we can solve the non-local Riemann Hilbert (NLRH) problem with scattering data $\vct{t}_R$ and $\rho$. This is a consequence of results of \cite{Novikov1992}. The NLRH problem can be solved with the formulas and equations of \cite[Theorem 6.1]{Novikov1992}. Since $|\vct{t}_R(\lambda)| \leq |\vct{t}(\lambda)|$ for every $\lambda \in \C$, the truncated scattering transform satisfies the estimates required in \cite[Theorems 6.1 and 6.2]{Novikov1992} (which are already satisfied by $\vct{t}$, thanks to our assumptions) in order to solve the NLRH problem via integral equations by iteration.\smallskip

Estimate \eqref{regular} will be a consequence of technical results in \cite{Santacesaria2015}, originally obtained to prove stability estimates for this problem. Since the NLRH problem can be solved for $q_0$ and $q_R$, we can repeat the arguments of \cite[Section 4]{Santacesaria2015} in order to obtain the equality \cite[identity (4.12)]{Santacesaria2015}
\begin{equation} \label{idest}
q_0(z)-q_R(z)=2 i \sqrt{E}(A-A_R +B-B_R + C-C_R),
\end{equation}
where $A, A_R, B,B_R,C,C_R$ are constructed as $A, B, C$ in \cite[Section 4]{Santacesaria2015} with $\mathrm{sgn}(|\lambda|^2-1)\frac{\vct{t}(\lambda)}{4\pi \bar \lambda}$ and $\mathrm{sgn}(|\lambda|^2-1)\frac{\vct{t}_R(\lambda)}{4\pi \bar \lambda}$ instead of $r(\lambda)$. In \cite[Section 5]{Santacesaria2015} the right hand side of \eqref{idest} is estimated in terms of scattering data. First, the estimate after \cite[estimate (5.6)]{Santacesaria2015} in the present notation reads
\begin{equation}
|A-A_R| \leq c(\Omega,N,m)\left(\sqrt{E}\left\| \left(\frac{1}{\lambda}+\bar \lambda\right)\frac{\vct{t}(\lambda)}{\bar \lambda}\right\|_{L^1(\C \setminus \A_R)}\! \! \! \! \! \!+\left\|\frac{\vct{t}(\lambda)}{\bar \lambda}\right\|_{L^p(\C \setminus \A_R)}\right),
\end{equation}
for some $p \in ]1,2[$, where $\A_R$ is the annulus defined in the statement. For $B-B_R$ we use the estimate after \cite[estimate (5.8)]{Santacesaria2015}. Since $q_0$ and $q_R$ correspond to scattering data $(\vct{t},\rho)$ and $(\vct{t}_R,\rho)$, the first term in this estimate vanishes; but also the third one - corresponding to $\delta r'_a$, defined in the statement of \cite[Proposition 4.2]{Santacesaria2015} - vanishes if we fix $a \leq 2$, since we chose $R >2$ (so $\vct{t} \equiv \vct{t}_R$ in the annulus $\A_2$). Thus we obtain the estimate
\begin{equation}
|B-B_R| \leq c(\Omega,N,m)\left\| \left(\frac{1}{\lambda}+\bar \lambda\right)\frac{\vct{t}(\lambda)}{\bar \lambda}\right\|_{L^{s,s'}(\C \setminus \A_R)},
\end{equation}
for some $1 < s<2<s'< +\infty$, where $\| \cdot \|_{L^{s,s'}} = \| \cdot \|_{L^s} + \| \cdot \|_{L^{s'}}$. The same argument applies to $C-C_R$ and we get
\begin{equation}\label{estc}
|C-C_R| \leq c(\Omega,N,m)\left\| \left(\frac{1}{\lambda}+\bar \lambda\right)\frac{\vct{t}(\lambda)}{\bar \lambda}\right\|_{L^{s,s'}(\C \setminus \A_R)}.
\end{equation}
Finally, we use \cite[Lemma 3.1]{Santacesaria2015}, which gives $L^p$ estimates of $\mathrm{sgn}(|\lambda|^2-1)\frac{\vct{t}(\lambda)}{4\pi \bar \lambda}$ (this corresponds to $r(\lambda)$ in that lemma) near $0$ and $\infty$. We have
\begin{equation}
\left\| |\lambda|^{j} \frac{\vct{t}(\lambda)}{\bar \lambda} \right\|_{L^p(\C \setminus \A_R)}\leq c(m,N)E^{-m/2}R^{-m+2},\quad \text{for } j =-1,0,1,
\end{equation}
for $R \geq R_0(N,m)$ and $p \geq 1$. This combined with \eqref{idest}-\eqref{estc} yields the main estimate \eqref{regular}.
\end{proof}

\begin{corollary} \label{corolregular}
Let $q$, $\Omega$, $m$, $N$, $R$, $E$ and $q_R$ be as in Lemma \ref{lemmaregular}. Let $\chi'_R$ be the characteristic function of the ball of radius $\sqrt{E}R$ centered in the origin and define $\tilde q_R = \F^{-1}[ \chi'_R \F q]$, where $\F$ is the 2D Fourier transform. Then there is a constant $C = C(\Omega,N,m) >0$ such that
\begin{equation}
\|q_R - \tilde q_R\|_{L^{\infty}(\Omega)} \leq C E^{-(m-2)/2}R^{-(m-2)}.
\end{equation}
\end{corollary}
\begin{proof}
Since $q_0 \in W^{m,1}(\Omega)$ we have $|\F q_0(w)| \leq c(\Omega,N,m)|w|^{-m}$ for $|w| \geq 1$. Then
\begin{align*}
|q_0(z)-\tilde q_R(z)|&= |\F^{-1}[(1-\chi'_R)\F q_0]|\\ 
&\leq c(\Omega,N,m)\int_{|w|\geq \sqrt{E}R} \frac{d\Ree w d\Imm w}{|w|^{m}} \leq \frac{c(\Omega,N,m)}{(\sqrt{E}R)^{m-2}},
\end{align*}
for every $z \in \Omega$. This combined with Lemma \ref{lemmaregular} gives the corollary.
\end{proof}


Now, choose an integer $N_\lambda>0$ and radii $1<R_1<R_2$. For spectral parameters $R_1<\abs{\lambda}<R_2$ define a $N_\lambda\times N_\lambda$ grid. For these values we precompute the matrices $\mathrm{S}_\lambda$ in order to solve the boundary integral equation. The radius $R_1>1$ is taking out values of $\lambda$ close to the unit circle, since we have problems in computing the Faddeev Green's function for these values. The use of $R_2$ was justified in the above Lemma and is analogous to the truncation radius of the zero-energy case acting as a regulation parameter. Depending on the case the value of $R_2$, outside of which the computational problems arise, changes. The computational problems can be seen from the computed scattering transform.

We denote by $\mathcal{F}^{-1}$ the transformation from the Fourier series domain to the function domain and simply use \eqref{scat-trans-DN} to get $\vct{t}(\lambda)$:
$$
\T(\lambda) = \int_{\bound}e^{\frac{i\sqrt{E}}{2}(\kk{\lambda}z+\kk{z}/\kk{\lambda})}\mathcal{F}^{-1}((\mathrm{L}_q-\mathrm{L}_{-E})\vct{\psi}_\lambda)ds.
$$
Recall the symmetry \eqref{tsymmetry} for any non-exceptional $\lambda$. Using this we can construct the scattering transform inside the unit circle. Depending on the scattering transform, choose the radius $R_2$ inside of which the numerical computation is usable. Then use the truncated scattering transform
\begin{equation}\label{truncation}
\vct{t}_R(\lambda) = \left\{
\begin{array}{rcl}
                      0,&\quad           &\abs{\lambda}\leq 1/R_2 \\
\vct{t}(1/\kk{\lambda}),& \quad 1/R_2\leq&\abs{\lambda}<1/R_1 \\
                      0,&\quad  1/R_1\leq&\abs{\lambda}\leq R_1 \\
       \vct{t}(\lambda),& \quad      R_1<&\abs{\lambda}< R_2\\
                      0,&\quad           &\abs{\lambda}\geq R_2.
\end{array}
\right.
\end{equation}

\begin{remark}
Although the above-defined truncated scattering transform differs from the one in Lemma \ref{lemmaregular}, a regularization estimate similar to \eqref{regular} (but less sharp) can be proved using the same ideas. Under the hypothesis of Lemma \ref{lemmaregular}, by \cite[Estimate (2.18c)]{Novikov1999} we have
$$ |\vct{t}(\lambda)| \leq C(N)(1+E(|\lambda|+|\lambda|^{-1})^2)^{-m/2}, \qquad \lambda \in \C,$$
which gives
\begin{equation*}
\left\| |\lambda|^{j} \frac{\vct{t}(\lambda)}{\bar \lambda} \right\|_{L^p(1/R_1 \leq |\lambda| \leq R_1)}  \! \! \! \! \! \! \leq O\left(E^{-m/2}(R_1-1)\right),\quad \text{for } j =-1,0,1,
\end{equation*}
where $p \geq 1$.
This, combined with the proof of Lemma \ref{lemmaregular}, yields a reconstruction error of the order $E^{-m/2} \max\left(O(R_1-1), O( R_2^{-(m-2)})\right)$.
\end{remark}

\subsection{Solving the D-bar  equation}
We can solve the periodic version of the integral equation \eqref{IEIS-compact}, without the term $\MM$, using $\vct{t}_R$ and the analog of the solver fully detailed in \cite{Knudsen2004}. 
We will deal with the following integral equation,
\begin{equation}\label{IEIS-modified}
\mu_R = 1-(\CC\TT_R)\mu_R,
\end{equation}
where $\TT_R$ is the operator of \eqref{T-def} with $\T_R(\lambda)$ instead of $\T(\lambda)$. 
Equation \eqref{IEIS-modified} is solved by periodization and using a matrix-free implementation of GMRES. See \cite[section 15.4]{Mueller2012} for details.


\subsection{Reconstructing the potential}

Let $z_r$ be the reconstruction point of our choosing. Let $dz$ be the finite difference and define the points $z_1 = z_r+dz$, $z_2 = z_r-dz$, $z_3 = z_r+i\cdot dz$, $z_4 = z_r-i\cdot dz$. Using the earlier described methods we can solve the corresponding CGO solutions $\mu_R^i = \mu_R(z_i,\lambda)$, $i=1,2,3,4$. We combine the equations \eqref{mu-asymptotics} and \eqref{q0-recon1}, omit the term $\OO(1/\lambda)$, use a finite $\lambda$ and finite difference method for the differentiation to get the approximate reconstruction equation
\begin{equation}\label{q0recon-discrete}
q_0(z_r) \approx \lambda\sqrt{E}\left(i\frac{\mu_R^1-\mu_R^2}{2dz}+\frac{\mu_R^3-\mu_R^4}{2dz}\right).
\end{equation} 
Note that the result is computed in a grid of parameters $\lambda$ (since \eqref{IEIS-modified} is). We compute an average of $q_0(z_r)$ over values corresponding to $\abs{\lambda}=R_2$.

\section{Numerical results}\label{NUMERICS}

\noindent
In Section \ref{sec:Green_validation} we test the algorithm for  $g_\lambda(z)$ by computing the CGO solutions $\mu(z,\lambda)$ with $\abs{\lambda}>1$ for some potentials. This is done by solving the Lippmann-Schwinger type equation \eqref{LS-lambda} using the numerical solution method described in \cite[section 14.3]{Mueller2012}. Evaluating numerically the D-bar equation \eqref{D-bar} allows us to assess the accuracy of the CGO solutions.

We compute in Section \ref{sec:exceptional}  the scattering transform for various radially symmetric potentials and observe the emergence of exceptional points. This is analogous to the zero-energy study \cite{Music2013}.

In Section \ref{sec:reconstruction} we test the full D-bar algorithm for reconstructing several test potentials from their approximate Dirichlet-to-Neumann maps.

In Section \ref{sec:comparison} we test our algorithm against the Novikov-Santacesaria algorithm of \cite{Novikov01012013} with different energies.

\subsection{Validation of the numerical Faddeev Green's function}\label{sec:Green_validation}

\subsubsection{Definition of potentials}
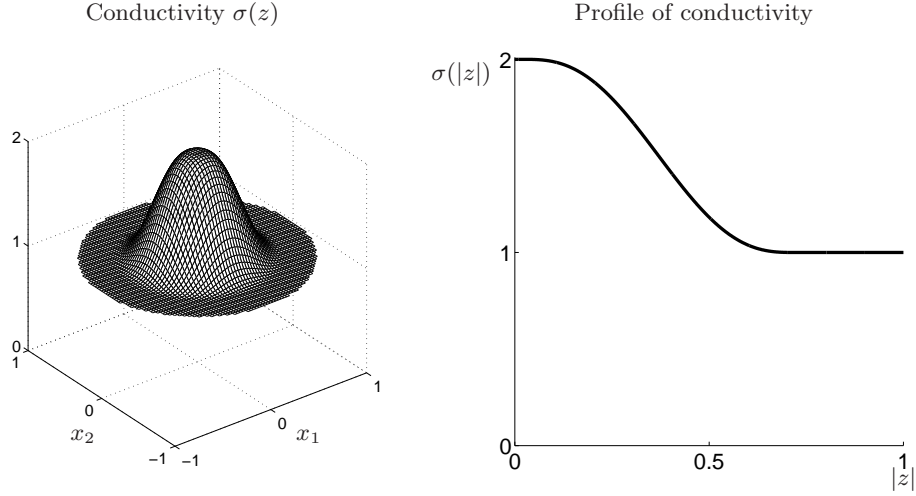
\begin{figure}
\unitlength=1mm
\begin{picture}(120,68)
 \put(0,0){\includegraphics[width=5cm]{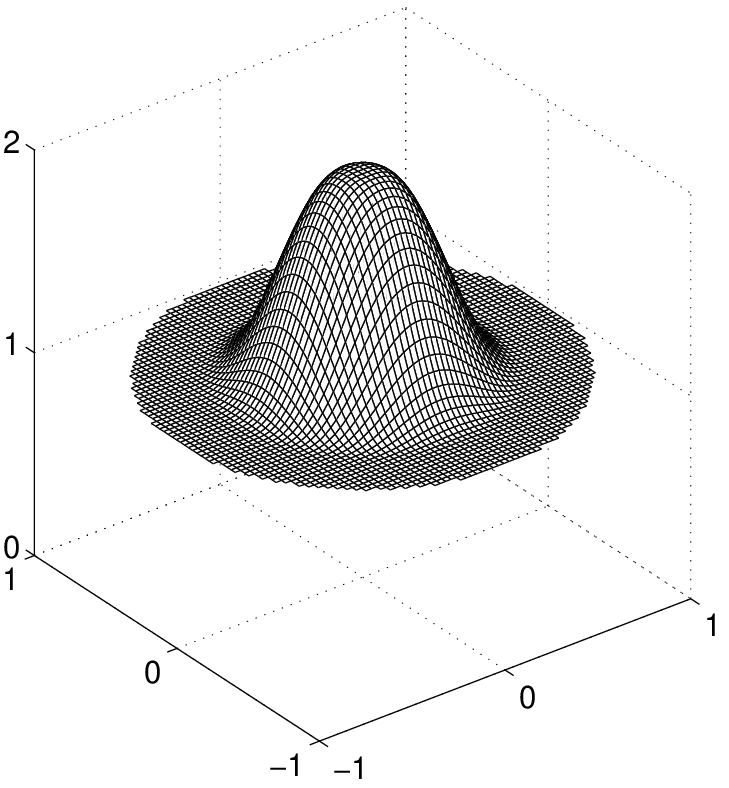}}
\put(10,60){\small Conductivity $\sigma(z)$}
 \put(65,0){\includegraphics[width=5.5cm]{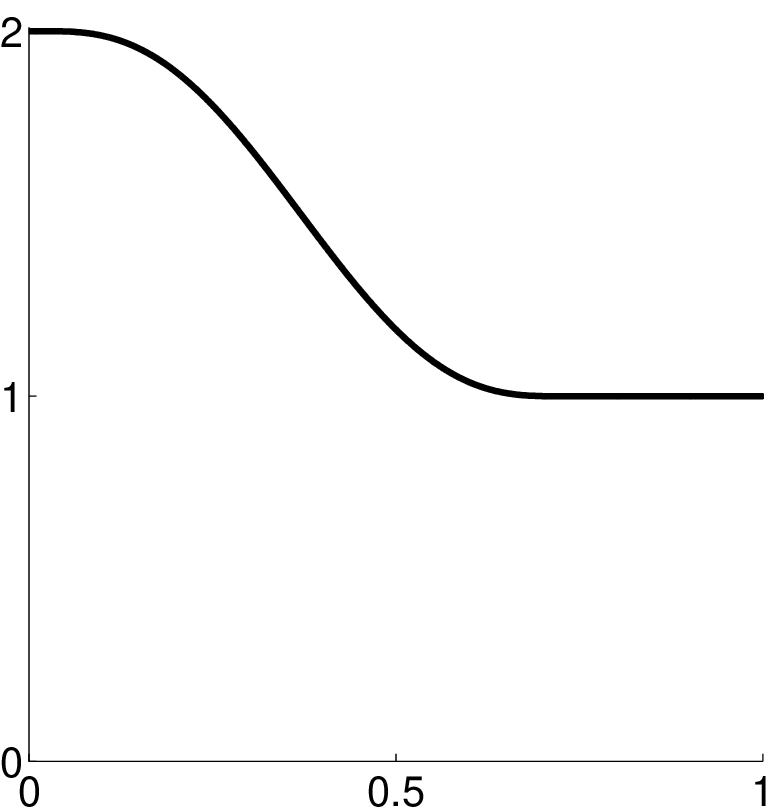}}
\put(75,60){\small Profile of conductivity }
\put(56,52){\small $\sigma(|z|)$}
\put(117,-2){\small $|z|$}
\put(38,4){\small $x_1$}
\put(8,4){\small $x_2$}
\end{picture}
\caption{\label{fig:sigma_plot}Mesh plot and profile plot of the rotationally symmetric conductivity $\sigma(z)=\sigma(|z|)$.}
\end{figure}

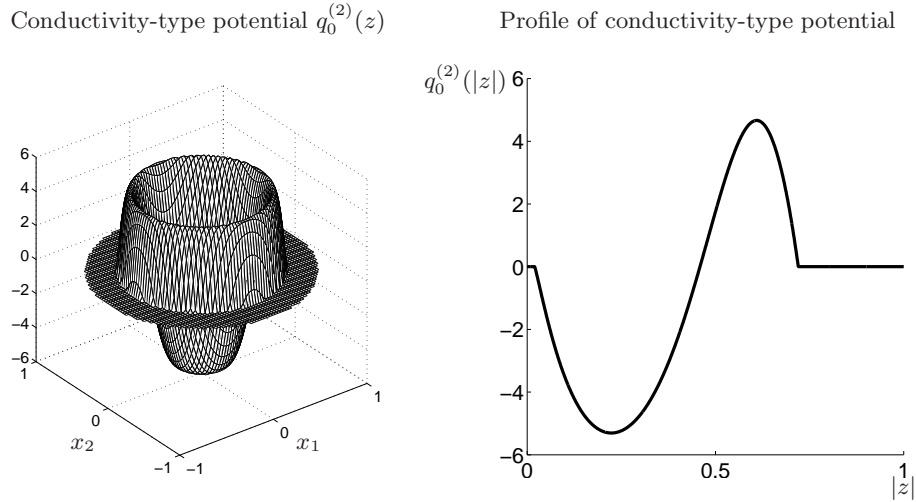
\begin{figure}
\unitlength=1mm
\begin{picture}(120,68)
 \put(0,0){\includegraphics[width=5cm]{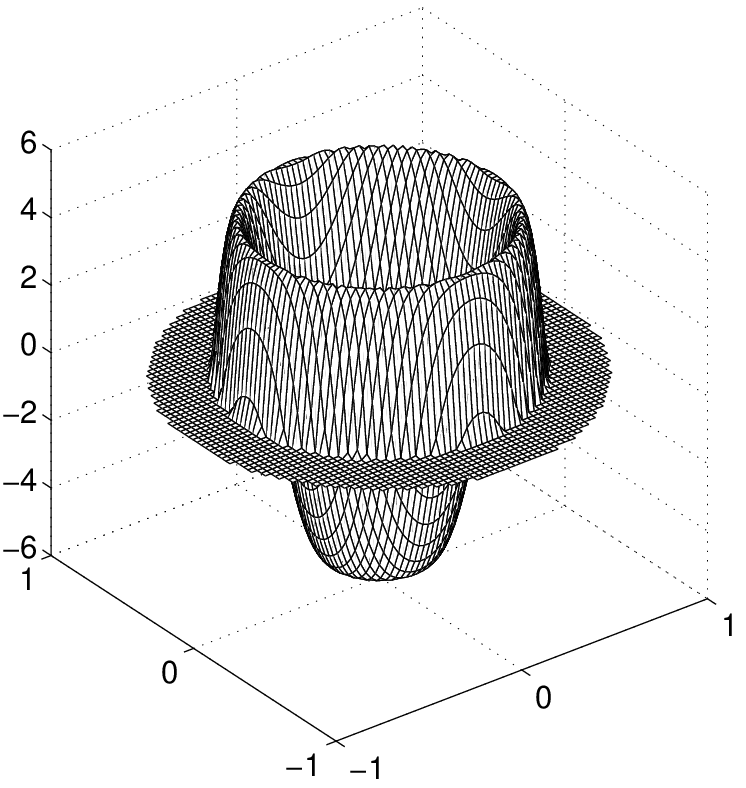}}
\put(0,60){\small Conductivity-type potential $q_0^{(2)}(z)$}
 \put(65,0){\includegraphics[width=5.5cm]{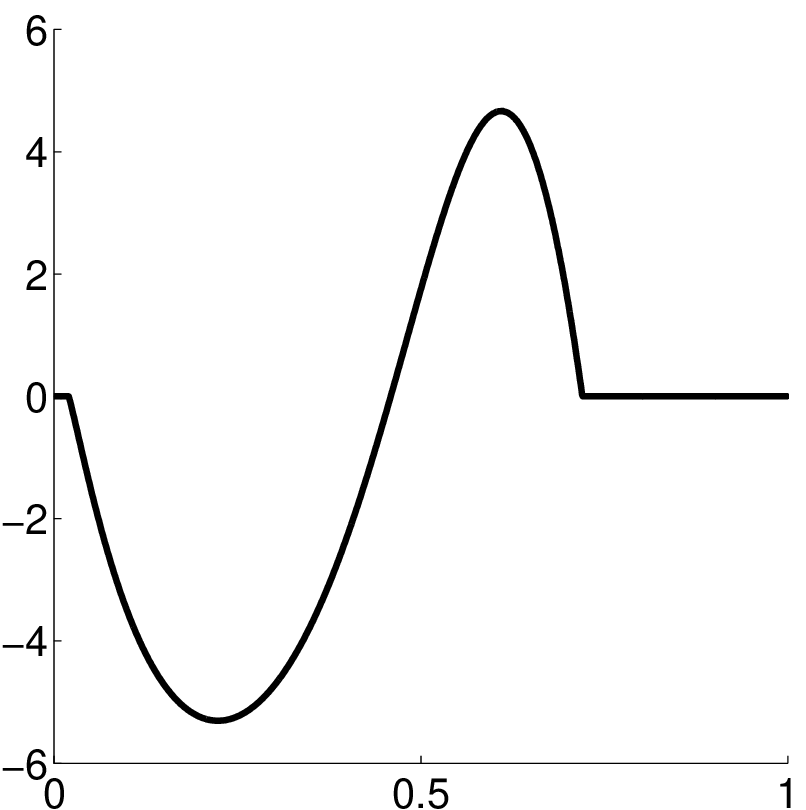}}
\put(65,60){\small Profile of conductivity-type potential }
\put(55,52){\small $q_0^{(2)}(|z|)$}
\put(117,-2){\small $|z|$}
\put(38,4){\small $x_1$}
\put(8,4){\small $x_2$}
\end{picture}
\caption{\label{fig:poten_plot}Mesh plot and profile plot of the conductivity-type potential $q_0^{(2)}(z)=q_0^{(2)}(\abs{z})$.}
\end{figure}

\begin{figure}
\unitlength=1mm
\begin{picture}(120,68)
 \put(0,0){\includegraphics[width=5cm]{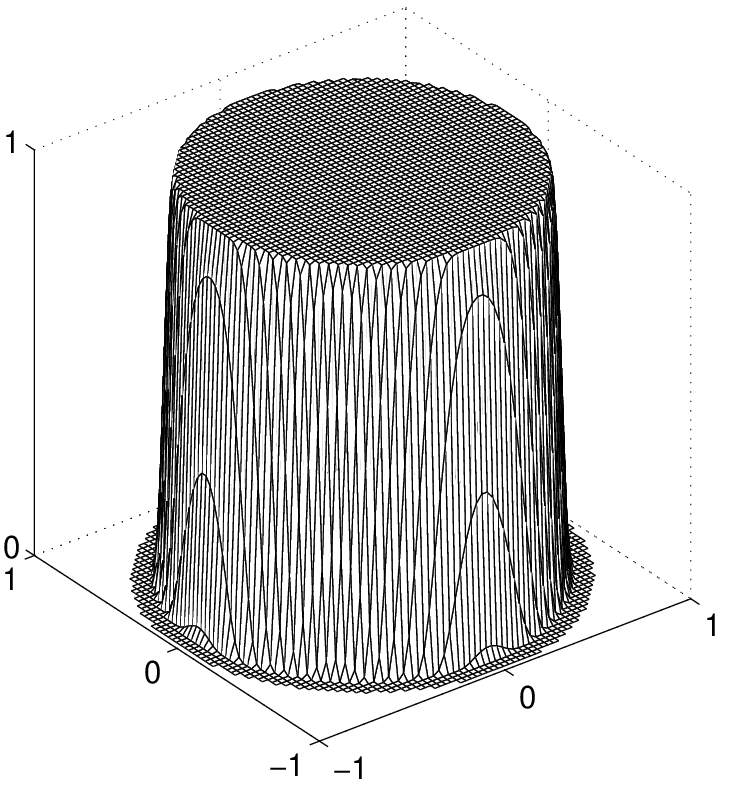}}
\put(10,60){\small Test function $\varphi(z)$}
 \put(65,0){\includegraphics[width=5.5cm]{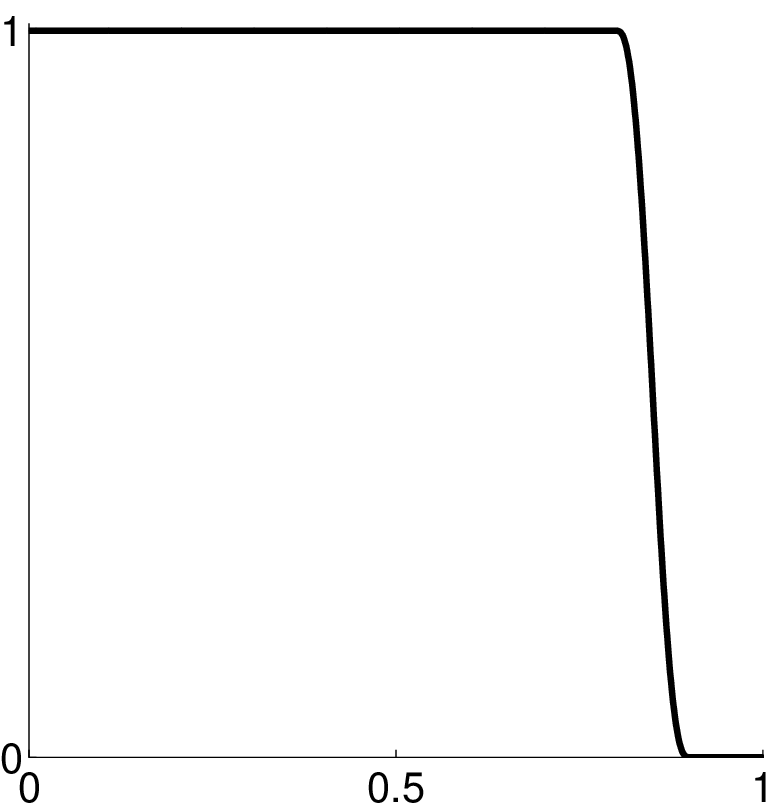}}
\put(75,60){\small Profile of test function }
\put(56,52){\small $\varphi(\abs{z})$}
\put(117,-2){\small $|z|$}
\put(38,4){\small $x_1$}
\put(8,4){\small $x_2$}
\end{picture}
\caption{\label{fig:testfunction_plot}Mesh plot and profile plot of the test function $\varphi(z)=\varphi(\abs{z})$.}
\end{figure}

We use exactly the same potentials as in the numerical part of \cite{Music2013}. Take radii $0<r_1<r_2<1$ and a polynomial $\tilde{p}(t) = 1-10t^3+15t^4-6t^5$. Set for $r_1\leq t \leq r_2$
\begin{equation*}
  p(t) = \tilde{p}(\frac{t-r_1}{r_2-r_1}).
\end{equation*}
Then, the approximate test function 
\begin{equation}
  \varphi(\abs{z}) = 
\left\{
\begin{array}{ccl}
1          & \mbox{ for } &0\leq \abs{z}\leq r_1 \\
p(\abs{z}) & \mbox{ for } &r_1< \abs{z}< r_2 \\
0          & \mbox{ for } &r_2\leq \abs{z}\leq 1,
\end{array}
\right.
\end{equation}
is in $C^2$. The values $r_1=0.8$ and $r_2=0.9$ were used. Consider  the radially symmetric potentials
\begin{eqnarray}
  q_{\alpha}^{(1)} &=& \alpha \varphi,\label{qalpha1-def}\\
  q_{\alpha}^{(2)} &=& \frac{\Delta\sqrt{\sigma}}{\sqrt{\sigma}}+\alpha\varphi, \label{qalpha2-def}
\end{eqnarray}
where $\alpha\in\R$ and $\sigma\in C^2(\Omega)$ with $\sigma\geq c>0$. The notion of "conductivity type potentials" is relevant also at positive energies as the problem of AT will include such a term in corresponding potential of the Gel'fand-Calder\'on problem. The mesh plot and the profile plot of the conductivity $\sigma$ are pictured in figure \ref{fig:sigma_plot}. The conductivity-type potential $q_0^{(2)}$ and the approximate test function $\varphi$ are pictured in figures \ref{fig:poten_plot} and \ref{fig:testfunction_plot} respectively.

\subsubsection{Verification of the computed CGO solutions}
\begin{figure}
\begin{picture}(120,60)
\epsfxsize=6cm
\put(1,0){\epsffile{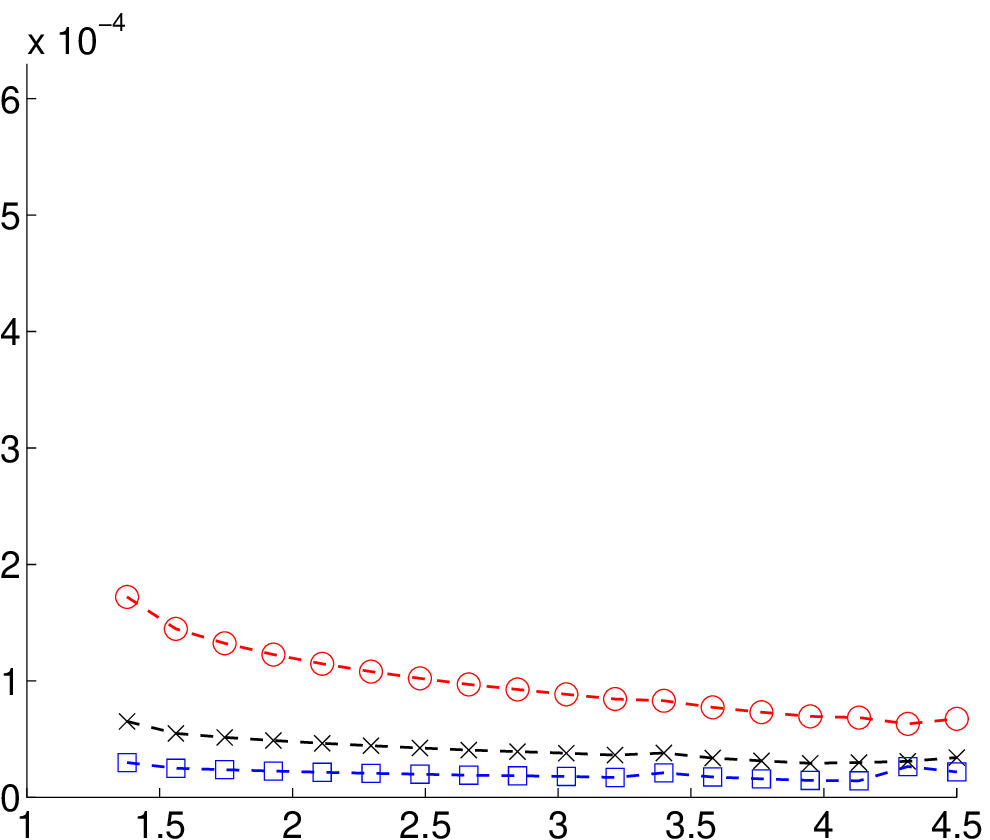}}
\epsfxsize=6cm
\put(63,0){\epsffile{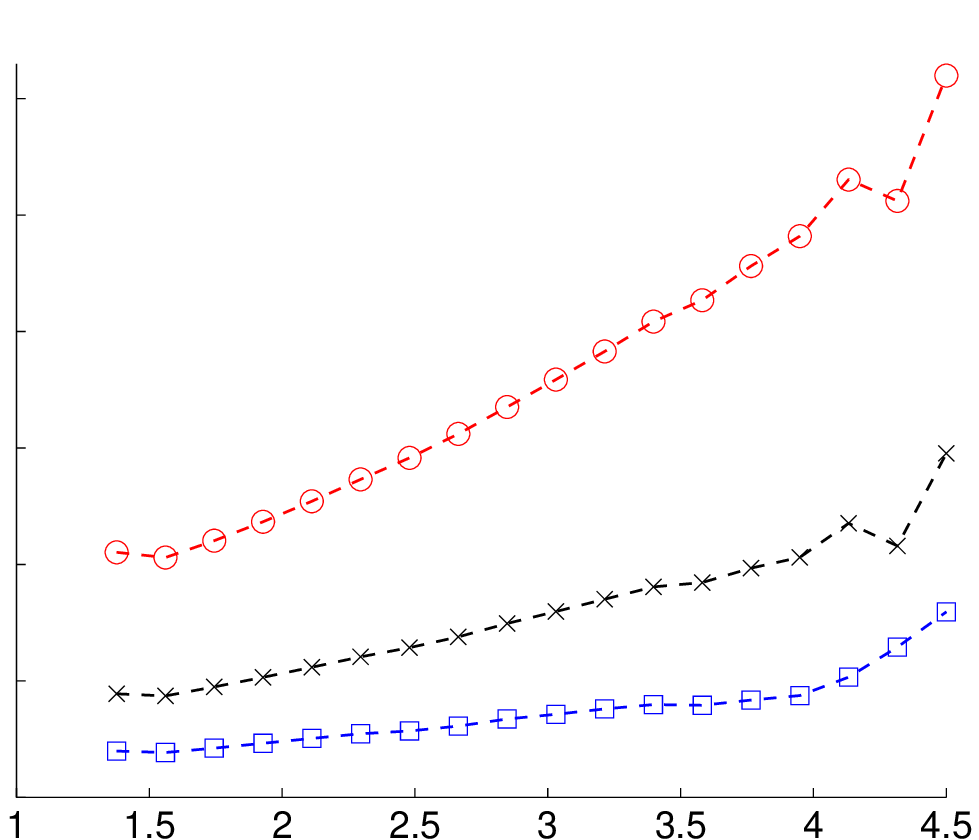}}
\put(123,-3){\small $\lambda$}
\put(0,54){\small $\norm{\cdot}_{L^2}$}
\end{picture}
\caption{\label{fig:dbar-error}Errors in the $\dbar$ -equation for two different potentials and different accuracies of the LS-solver; horizontal axis is the spectral parameter $\lambda$, vertical axis is the norm $\norm{\dbar\mu-\frac{1}{4\pi\kk{\lambda}}\vct{t}(\lambda)e_{-\lambda}(z)\kk{\mu_0}}_{L^2}$. On the left we used $q_0^{(2)}$ and on the right $q_{35}^{(2)}$. In red using circles is $M=7$, in black using crosses $M=8$ and in blue using squares is $M=9$. Two smallest values for $\lambda$ were omitted, for $\lambda=1.01$ the magnitude of the error was between 3 and 13, for the second smallest it was between 0.003 and 0.02.}
\end{figure}

In this section we fix $E=1$. To verify that the CGO solutions (and subsequently the Faddeev Green's function) are correct, we test the $\dbar$ -equation \eqref{D-bar} using the five-point stencil method with the finite difference $d\lambda=0.0001$. Take parameters $\lambda$ from 1.01 to 4.5. Take the potentials $q_0^{(2)}$ and $q_{35}^{(2)}$ to test two very differently sized potentials. For each $\lambda=\lambda_1+\lambda_2 i$, compute
\begin{enumerate}
\item The CGO solution $\mu_0$ in the $z$-grid, corresponding to the parameter $\lambda$.
\item The CGO solutions $\mu_1,\mu_2,\mu_3,\mu_4,\mu_5,\mu_6,\mu_7$ and $\mu_8$ using $\lambda+d\lambda$, $\lambda+2d\lambda$, $\lambda-d\lambda$, $\lambda-2d\lambda$, $\lambda+d\lambda i$, $\lambda+2d\lambda i$, $\lambda-d\lambda i$ and $\lambda-2d\lambda i$ respectively.
\item The functions $e_\lambda(z)$ and $e_{-\lambda}(z)$.
\item The scattering transform $\vct{t}(\lambda)$ of \eqref{scattering-transform} with $\mu = \mu_0$, $q_0 = q_0^{(2)}$ and $q_0 = q_{35}^{(2)}$.
\item The derivatives and the $\dbar$ -operation by
\begin{eqnarray*}
\partial_{\lambda_1}\mu &=& \frac{-\mu_2+8\mu_1-8\mu_3+\mu_4}{12d\lambda}\\
\partial_{\lambda_2}\mu &=& \frac{-\mu_6+8\mu_5-8\mu_7+\mu_8}{12d\lambda}\\
\dbar \mu &=& \frac{1}{2}(\partial_{\lambda_1}+i\partial_{\lambda_2})\mu.
\end{eqnarray*}
\item The error
  \begin{equation}\label{dbar-error}
    \norm{\dbar \mu-\frac{1}{4\pi\kk{\lambda}}\vct{t}(\lambda)e_{-\lambda}(z)\kk{\mu_0}}_{L^2(D(0,1))}.
  \end{equation}
\end{enumerate}
The above computations of CGO solutions are done with the solution algorithm described in \cite[section 14.3]{Mueller2012}; it is called {\em LS-solver} below. The $z$-grid has $2^M{\times}2^M$ points. In figure \ref{fig:dbar-error} we see the error \eqref{dbar-error} as a function of $\lambda$ using $q_0^{(2)}$ on the left, $q_{35}^{(2)}$ on the right. The parameter $M$ is increased from 7 to 9. As expected, the error decreases as $M$ increases as it increases the accuracy of the LS-solver. The smallest values of $\lambda$ were omitted in the pictures, for $\lambda=1.01$ the magnitude of the error was between 3 and 13, for the second smallest $\lambda$ it was between 0.003 and 0.02. For values of $\lambda$ near $\abs{\lambda}=1$ the numerical method of $g_\lambda(z)$ has great error due to very small value of $k_2$.

In the reconstruction of the potential we use values as large as $\abs{\lambda}\approx 15$. Not pictured here, this test was done also for larger values of $\lambda$, the error seems to be of similar magnitude for any $1.01<\abs{\lambda}<15$. In conclusion, the method for computing $g_\lambda(z)$ is valid and accurate enough for our purposes.

\subsection{Numerical investigation of exceptional points}\label{sec:exceptional}
In this section we fix $E=1$. For a given potential there may be values of parameter $\lambda$ for which there exists no unique CGO solution. Such $\lambda$ values are called exceptional points. We follow here the zero-energy study \cite{Music2013} and compute numerically CGO solutions at posivite energy. Exceptional points will show up as singularities in computation.

\begin{figure}
\begin{picture}(120,65)
\epsfxsize=6cm
\put(0,5){\epsffile{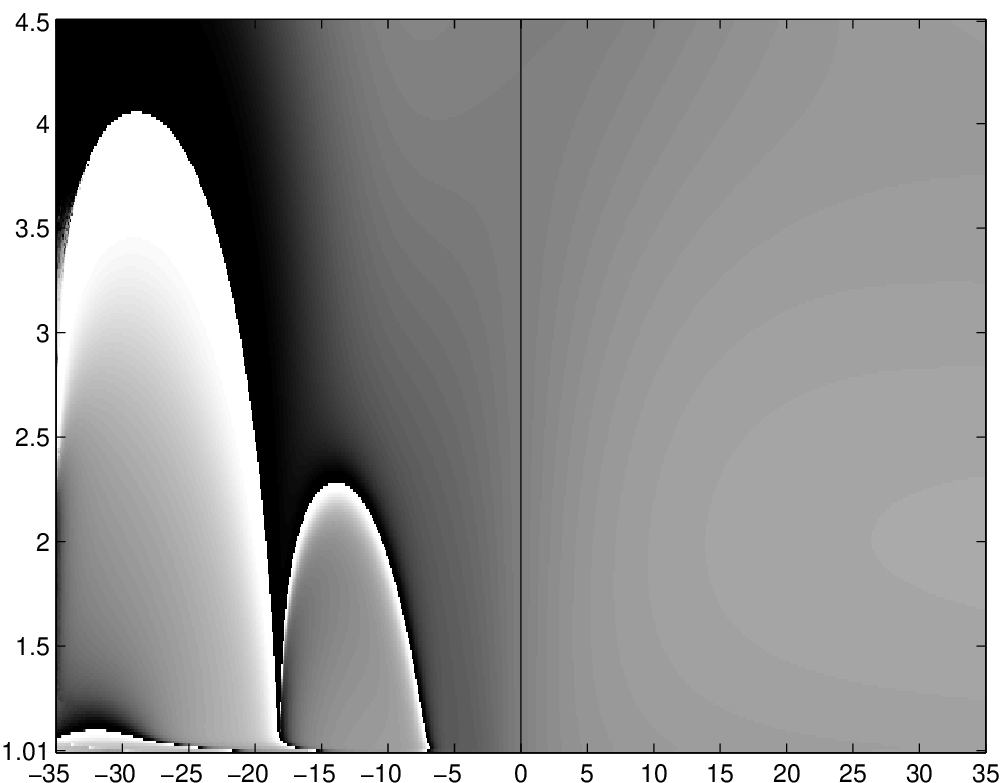}}
\epsfxsize=6cm
\put(60,5){\epsffile{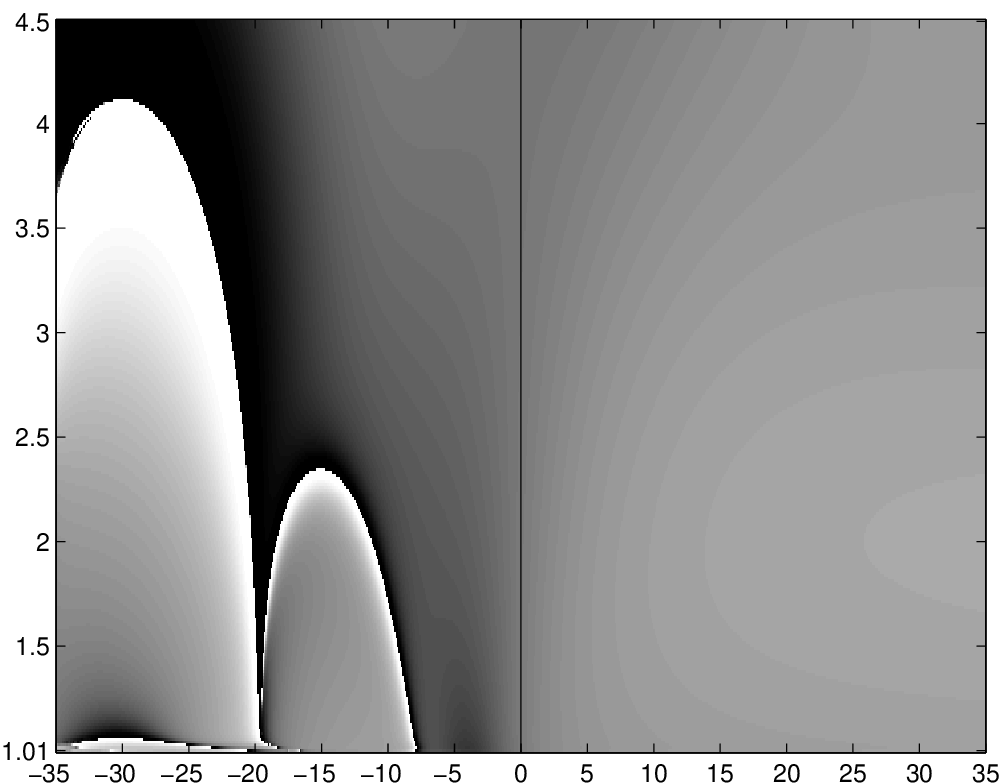}}
\put(115,2){\small $\alpha$}
\put(0,53){\small $\abs{\lambda}$}
\put(20,54){$\vct{t}(\abs{\lambda})$ for $q_\alpha^{(1)}$}
\put(80,54){$\vct{t}(\abs{\lambda})$ for $q_\alpha^{(2)}$}
\end{picture}
\caption{Scattering transform for the potential $q_\alpha^{(1)} = \alpha\varphi$ on the left, for the potential $q_\alpha^{(2)} =\Delta\sqrt{\sigma}/\sqrt{\sigma}+\alpha \varphi$ on the right. The x-axis is $\alpha = -35\ldots35$, y-axis is $\lambda = 1.01\ldots4.5$. Compare to figures 3 and 9 in \cite{Music2013}.\label{q-plane}}
\end{figure}
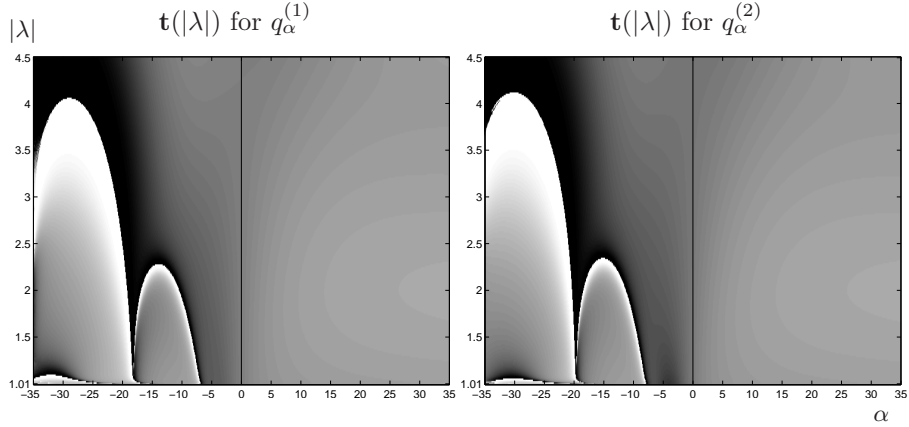

Recall the rotationally symmetric potentials $q_\alpha^{(1)}$ and $q_\alpha^{(2)}$ from \eqref{qalpha1-def} and \eqref{qalpha2-def}. See figures \ref{fig:sigma_plot},\ref{fig:poten_plot} and \ref{fig:testfunction_plot}. We use 250 discrete points of $\lambda$ and 701 discrete points of $\alpha$,
$$
\lambda = 1.01,\ldots,4.5,\quad \alpha = -35,\ldots,35.
$$
We use $M=8$ for the LS-solver (see \cite[section 14.3]{Mueller2012}) leading to $2^M\times 2^M$ sized $z$-grid. In figure \ref{q-plane} we plot the radially symmetric and real-valued scattering transform $\vct{t}(\lambda)=\vct{t}(|\lambda|)$ for the potential $q_\alpha^{(1)} = \alpha \varphi$ on the left, for the potential $q_\alpha^{(2)} =\Delta\sqrt{\sigma}/\sqrt{\sigma}+\alpha \varphi$ on the right. The $x$-axis is the parameter $\alpha$ and the $y$-axis is the modulus $|\lambda|$ of the spectral parameter. Black color represents very small negative values, and white very large positive values of $\vct{t}(\lambda)$. The lines where it abruptly changes between these colors are exceptional circles that move as the parameter $\alpha$ changes.

\begin{figure}
\unitlength=1mm
\begin{picture}(120,130)
 \put(10,88){\includegraphics[width=6cm]{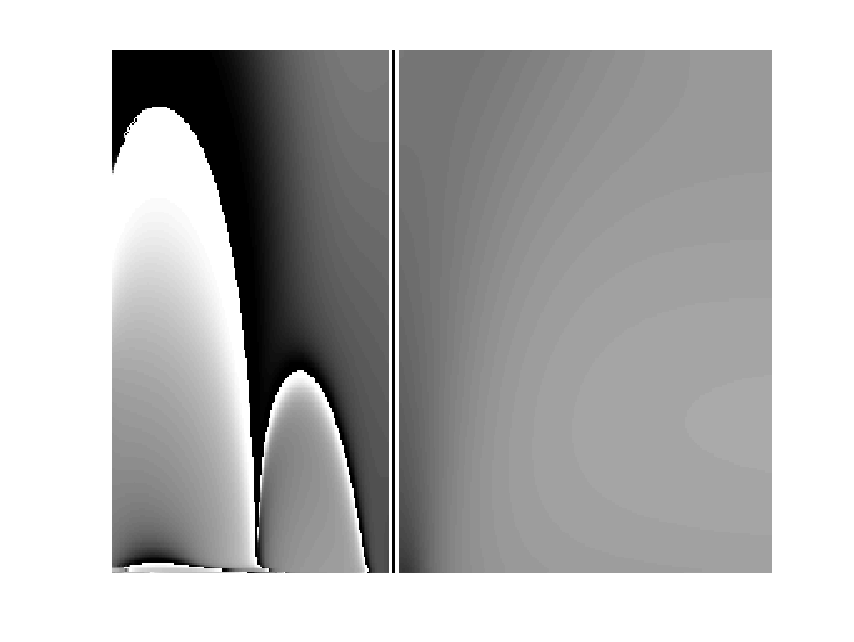}}
 \put(73,90){\includegraphics[width=5cm]{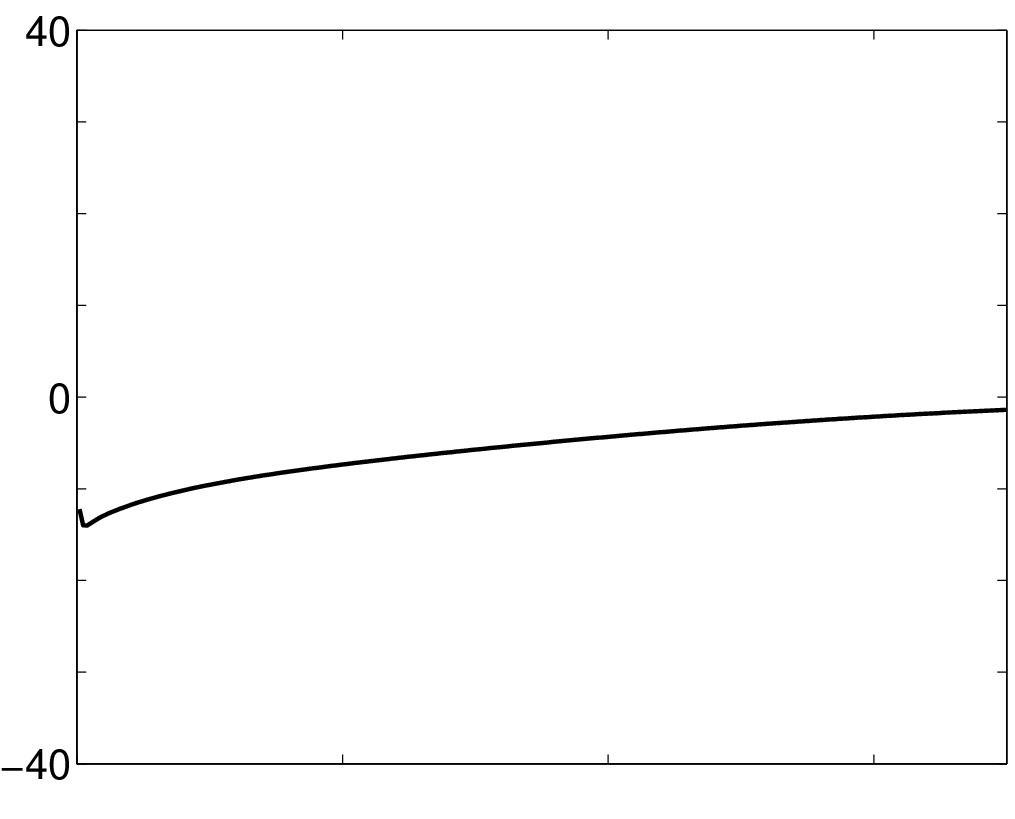}}
\put(0,110){$\alpha=-5$}
 \put(10,43){\includegraphics[width=6cm]{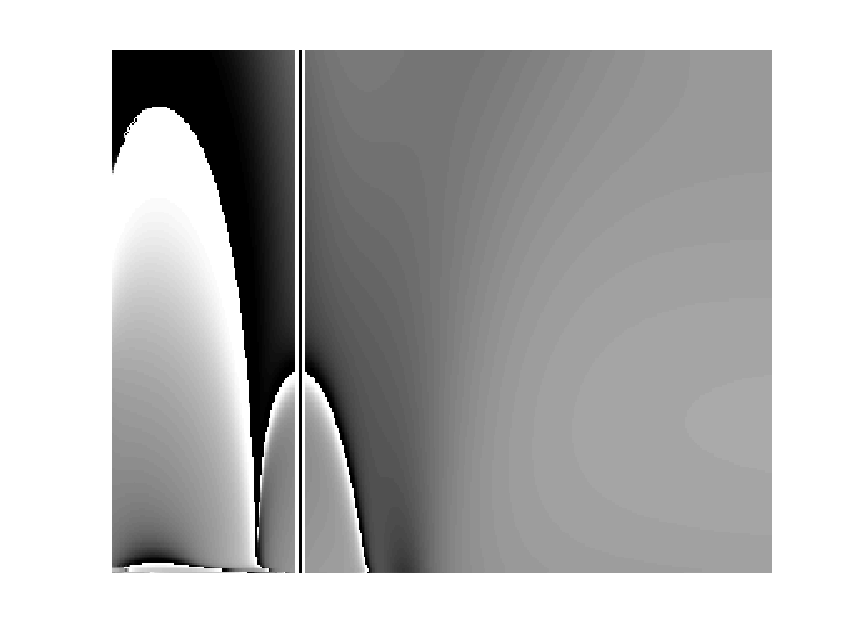}}
 \put(73,45){\includegraphics[width=5cm]{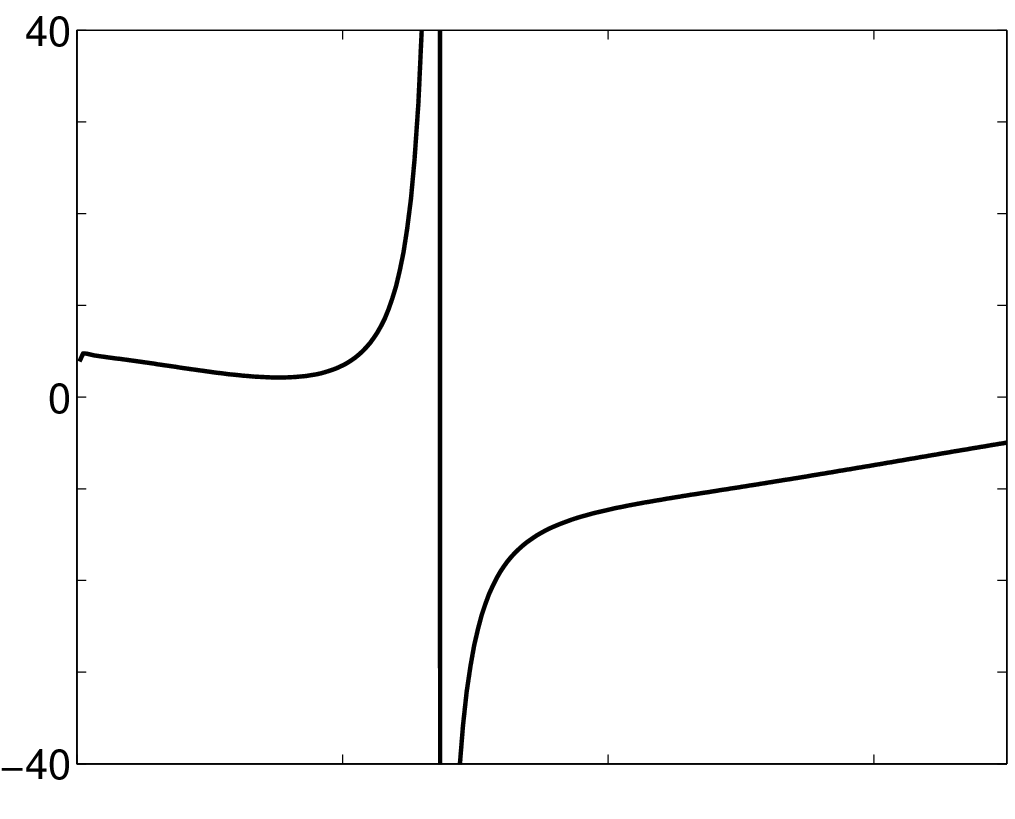}}
\put(0,65){$\alpha=-15$}
 \put(10,-2){\includegraphics[width=6cm]{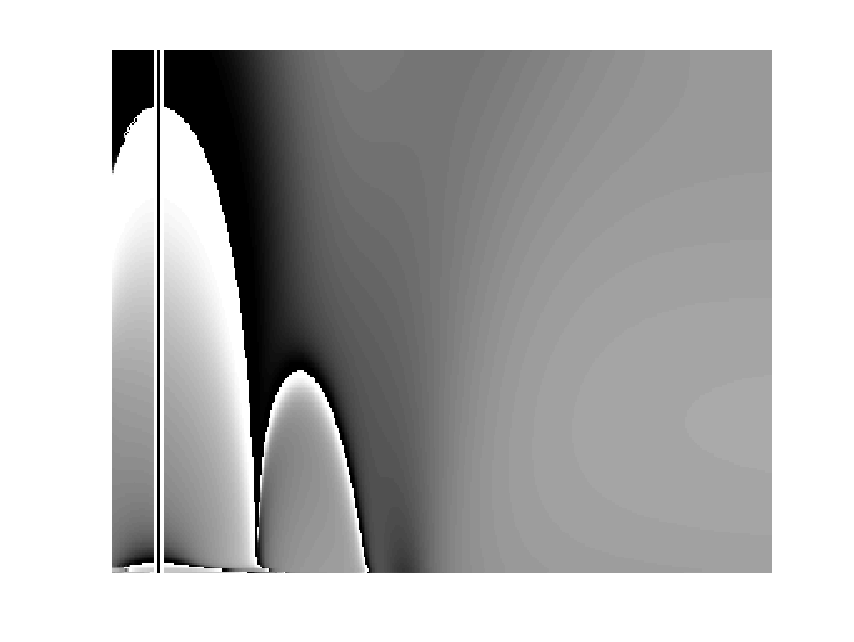}}
 \put(73,0){\includegraphics[width=5cm]{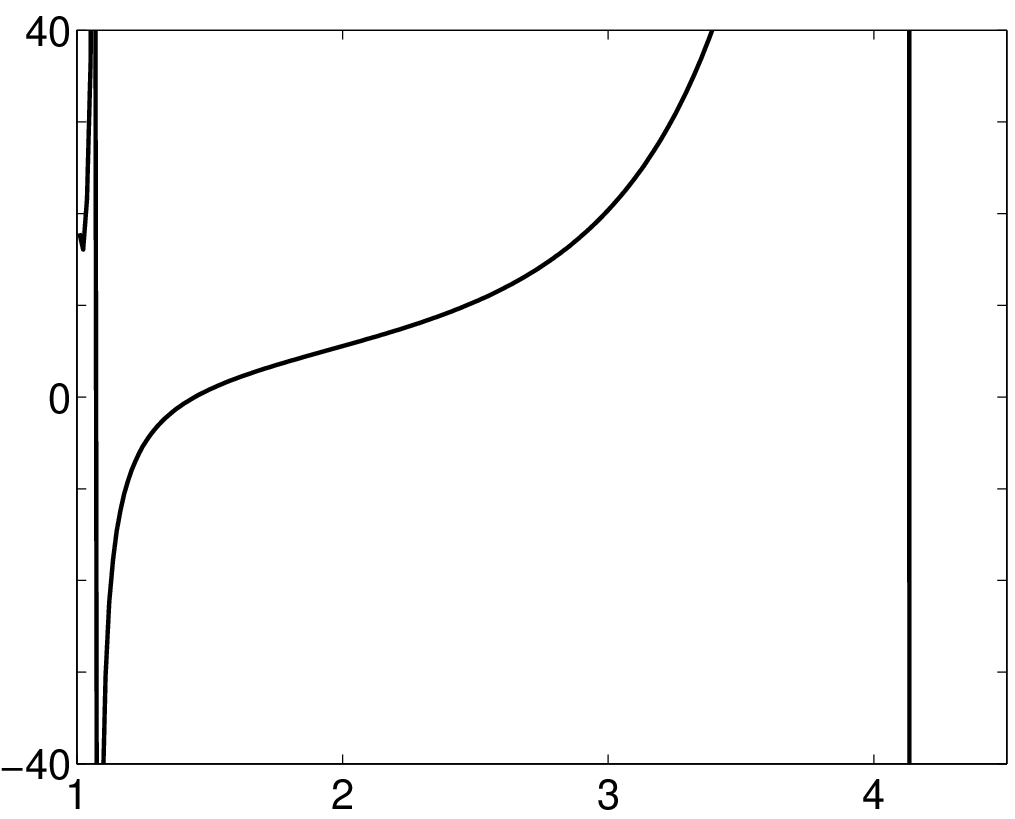}}
\put(0,20){$\alpha=-30$}
\put(117,-2){$\abs{\lambda}$}
\end{picture}
\caption{\label{scat_samples}On the right: the profile of $\vct{t}(\lambda)$ as a function of $\lambda$ using the potential $q_\alpha^{(1)} = \alpha\varphi$ with three different values of $\alpha$. On the left: the plane $\vct{t}(\lambda)$ for all parameters $\alpha$ with an indication of the location of the profile on the right. Compare to figure 4 in \cite{Music2013}.}
\end{figure}

In figure \ref{scat_samples} we plot the profile of the scattering transform $\vct{t}(\lambda)$ as a function of $\lambda$, using the potential $q_\alpha^{(1)} = \alpha \varphi$, with the values $\alpha=-5,-15,-30$. The exceptional circles can be seen as singularities in the profiles.

\subsection{Reconstructions of $q_0$}\label{sec:reconstruction}
\begin{figure}
\unitlength=1mm
\begin{picture}(120,65)
 \put(0,0){\includegraphics[width=6cm]{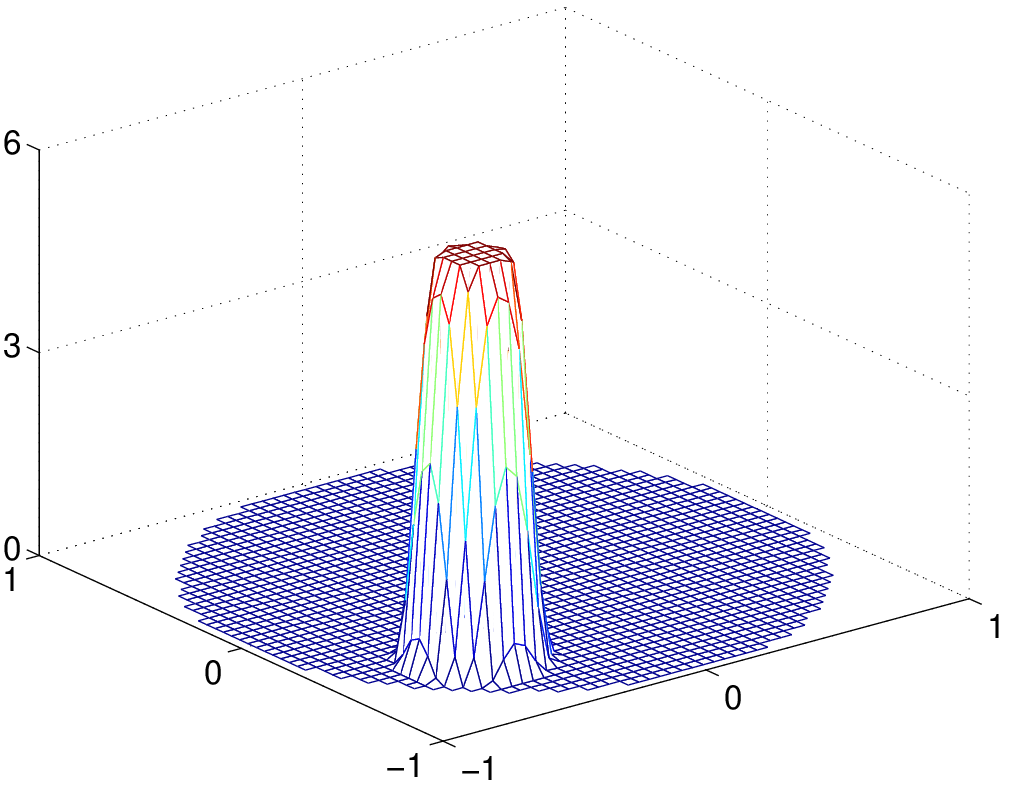}}
 \put(65,0){\includegraphics[width=5.5cm]{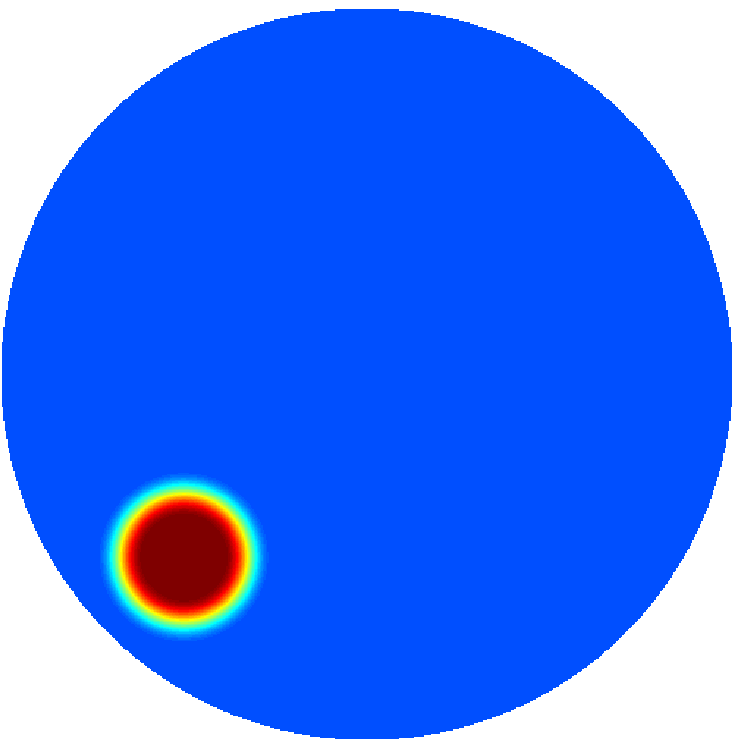}}
\put(40,2){\small $x_1$}
\put(8,4){\small $x_2$}
\end{picture}
\caption{\label{fig:q1plot}Mesh plot and 2D plot of the Case 1 potential $q_0(z)$.}
\end{figure}

\begin{figure}
\unitlength=1mm
\begin{picture}(120,65)
 \put(0,0){\includegraphics[width=6cm]{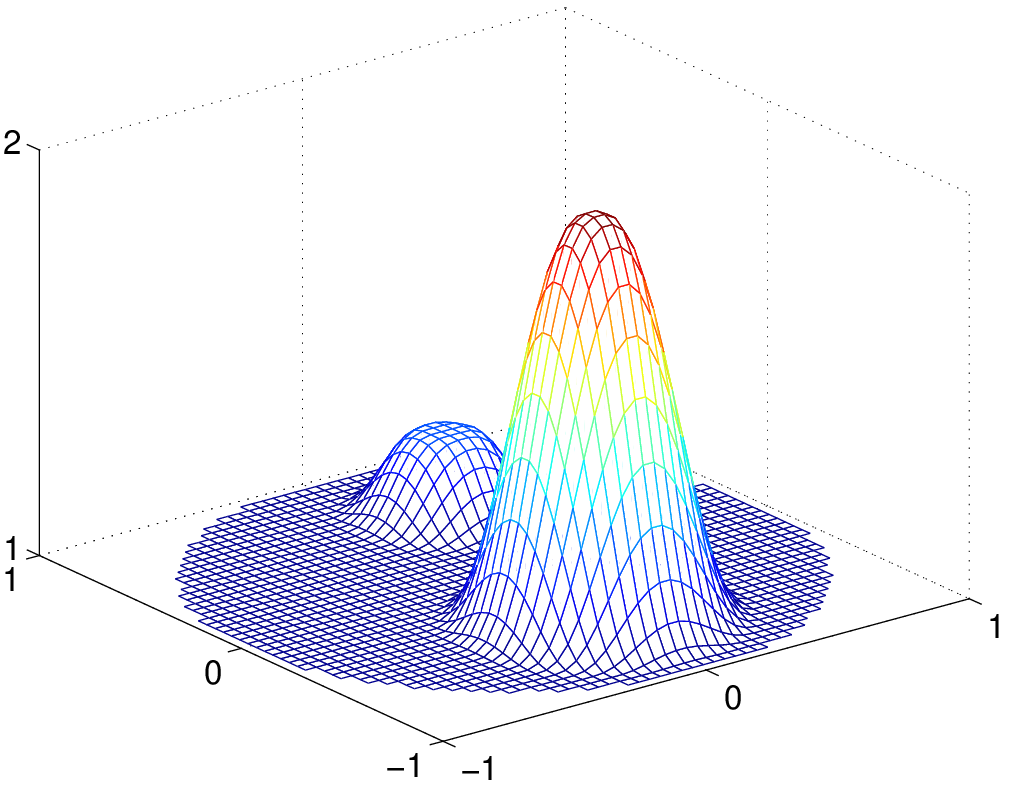}}
 \put(65,0){\includegraphics[width=5.5cm]{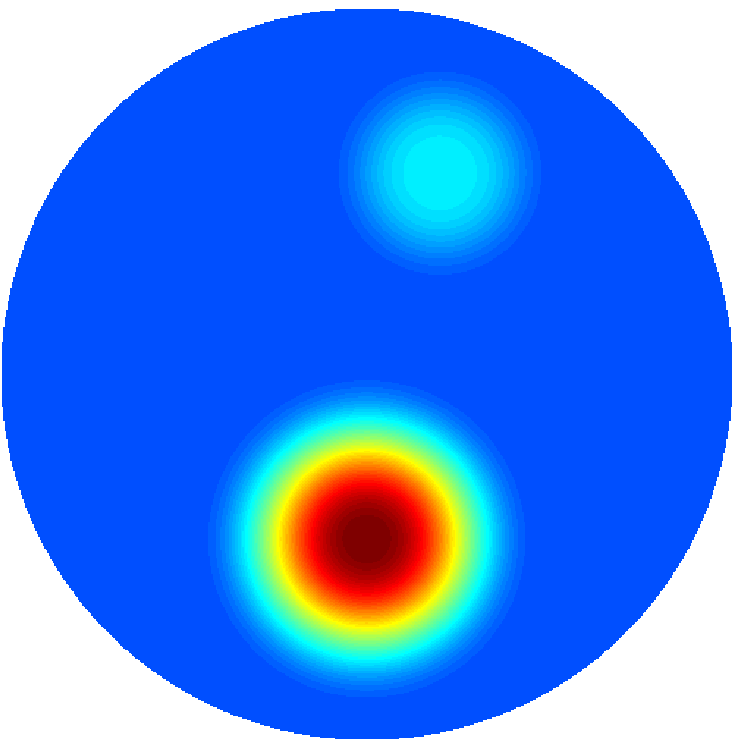}}
\put(40,2){\small $x_1$}
\put(8,4){\small $x_2$}
\end{picture}
\caption{\label{fig:sigmaplot}Mesh plot and 2D plot of the Case 2 conductivity $\sigma(z)$.}
\end{figure}

\begin{figure}
\unitlength=1mm
\begin{picture}(120,65)
 \put(0,0){\includegraphics[width=6cm]{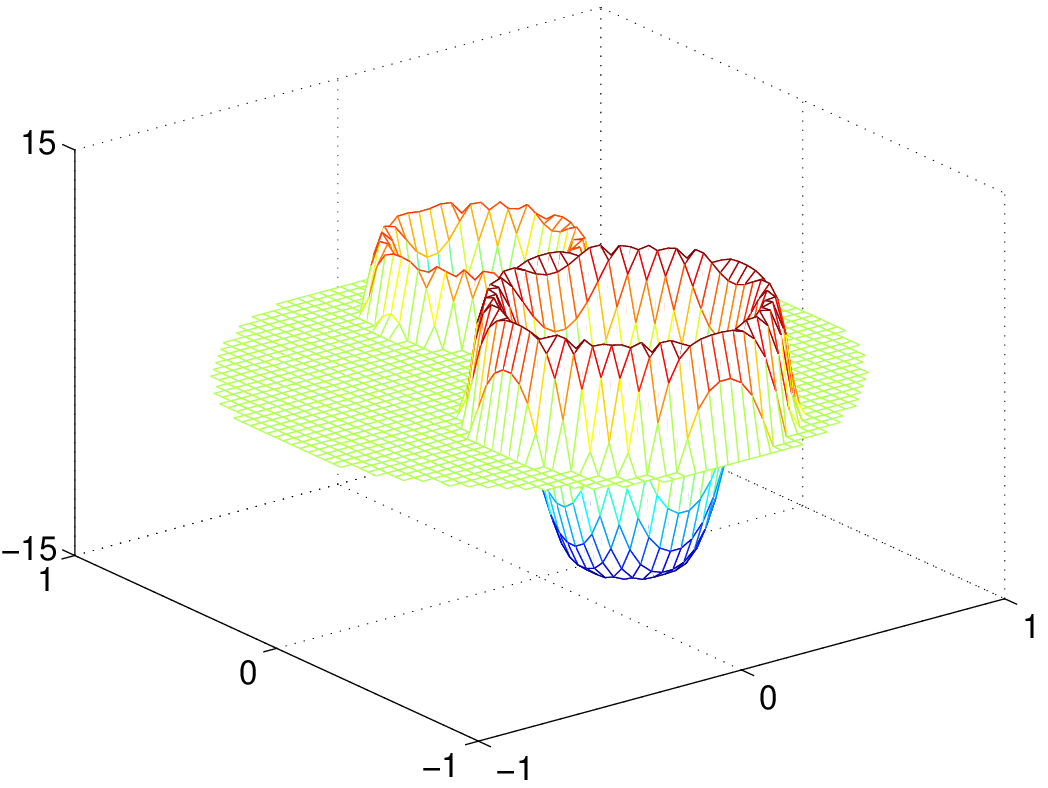}}
 \put(65,0){\includegraphics[width=5.5cm]{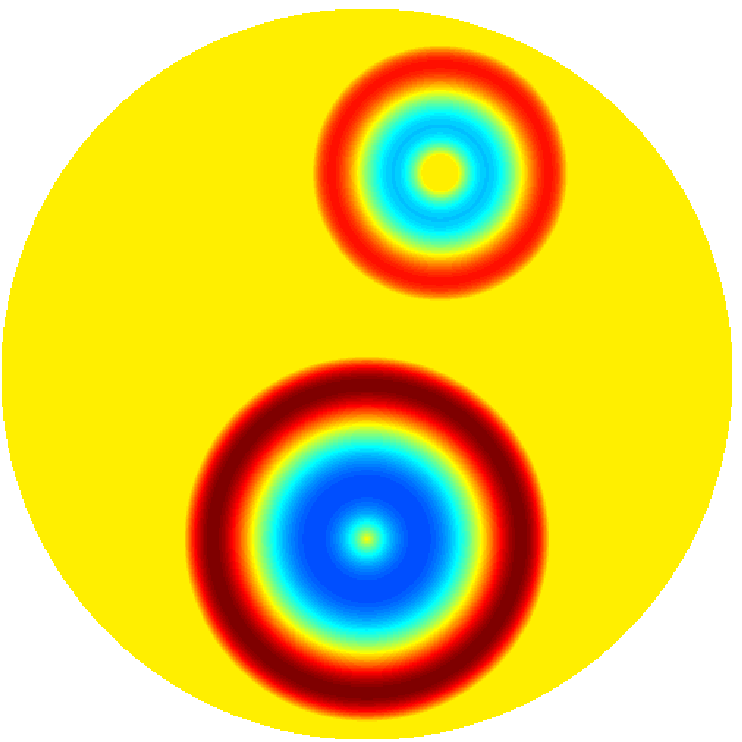}}
\put(40,2){\small $x_1$}
\put(8,4){\small $x_2$}
\end{picture}
\caption{\label{fig:q3plot}Mesh plot and 2D plot of the Case 2 potential $q_0(z) = \sigma(z)^{-1/2}\Delta\sigma(z)^{1/2}$.}
\end{figure}

In this section we fix $E=10^{-3}$. We reconstruct two non-symmetric potentials. The first is pictured in figure \ref{fig:q1plot}. The second is of conductivity type, the conductivity $\sigma$ is pictured in figure \ref{fig:sigmaplot} and the corresponding potential $q_0 = \sigma^{-1/2}\Delta\sigma^{1/2}$ in figure \ref{fig:q3plot}.

\subsubsection{Choice of parameters}
For the DN- and $\mathrm{S}_\lambda$ -matrices we use $N=16$, see \eqref{Lq-definition} and \eqref{Slambda-matrix}. We add gaussian noise to each element with \eqref{noiselevel} so that the relative matrix norm between the original DN-matrix and the noisy DN-matrix is 0.005\%. In the mesh for the FEM we have 1048576 triangles.

\begin{figure}[h]
\begin{picture}(120,128)
\epsfxsize=5.5cm
\put(0,65){\epsffile{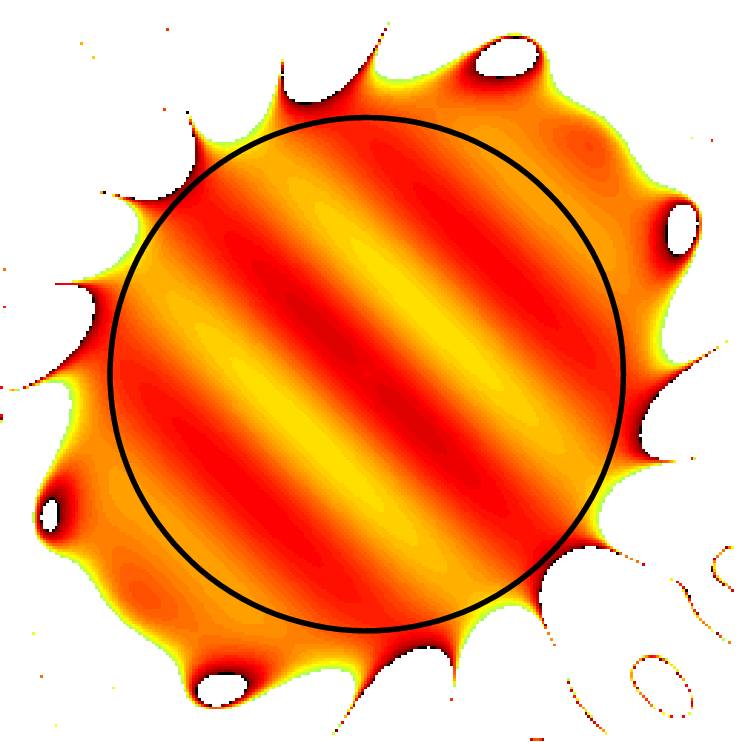}}
\epsfxsize=5.5cm
\put(60,65){\epsffile{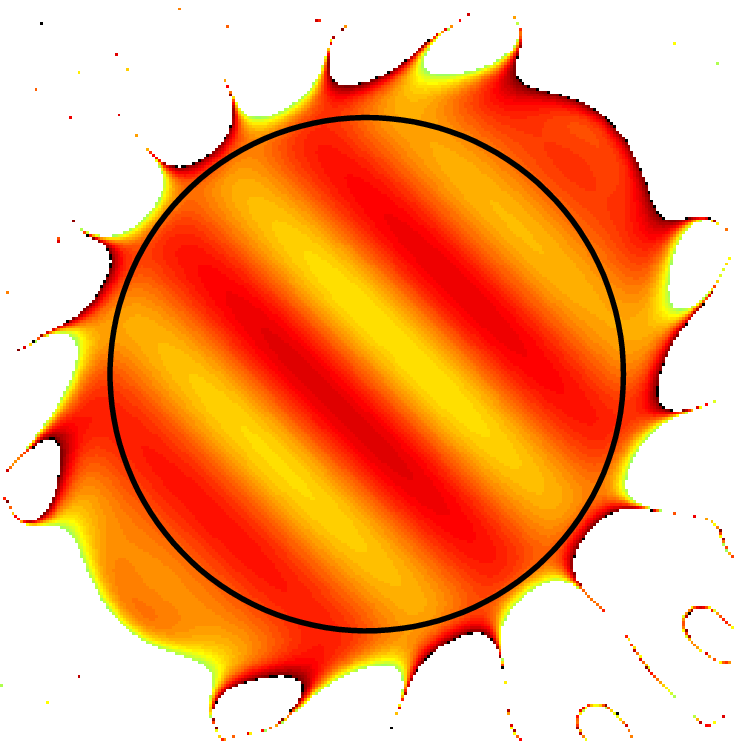}}
\epsfxsize=5.5cm
\put(0,0){\epsffile{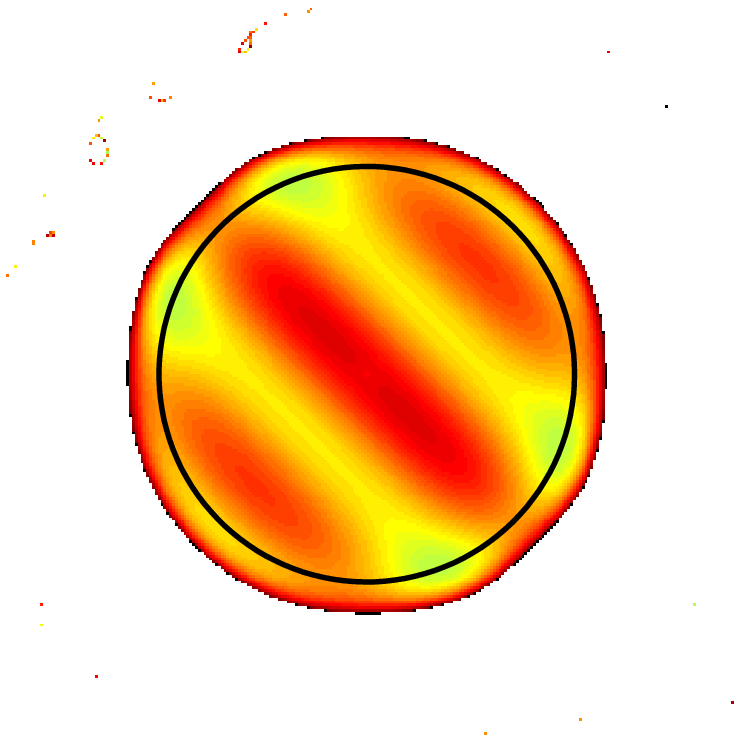}}
\epsfxsize=5.5cm
\put(60,0){\epsffile{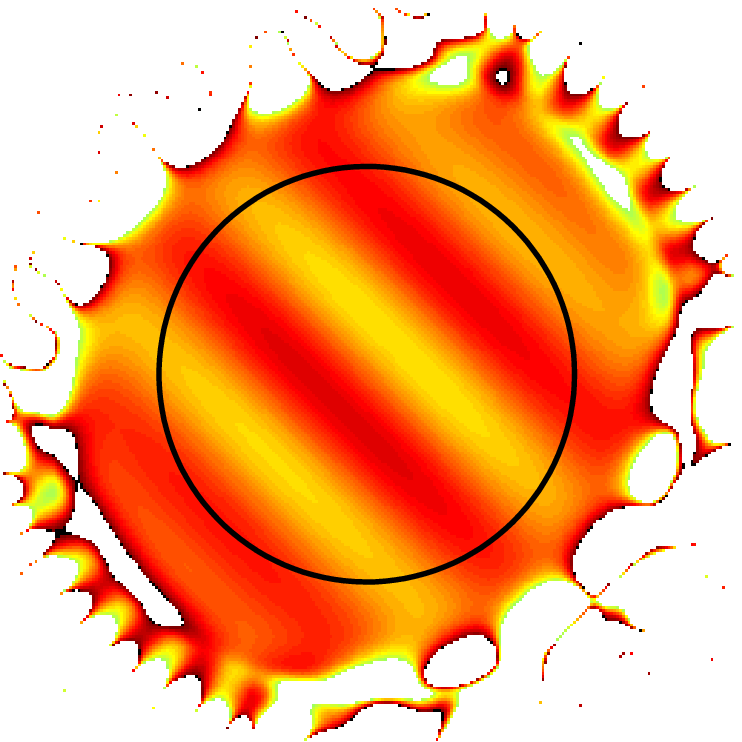}}
\put(21,124){$\Ree(\vct{t}(\lambda))$}
\put(85,124){$\Imm(\vct{t}(\lambda))$}
\put(48,124){Non-noisy, $\mathrm{L}_q$}
\put(48,48){Noisy, $\mathrm{L}^{\epsilon}_q$}
\end{picture}
\caption{\label{scats2D_DN}The scattering transform $\vct{t}(\lambda)$ of Case 1, the non-symmetric potential of figure \ref{fig:q1plot}. Real part on the left, imaginary part on the right, in a $\lambda$ -grid $[-600,600]\times [-600,600]i$. On the top row: the non-noisy DN-matrix $\mathrm{L}_q$ was used. On the bottom row: the noisy DN-matrix $\mathrm{L}^{\epsilon}_q$ was used. In the white areas the computation breaks down. The black line indicates the largest usable circle for the truncation $\vct{t}_R(\lambda)$. Energy level $E=0.001$.}
\end{figure}

Depending on the case, we cut off non-usable parts of the scattering transform. As an example, in figure \ref{scats2D_DN} we have the real and imaginary parts of the scattering transform of Case 1 potential computed using the non-noisy DN-matrix $\mathrm{L}_q$ and the noisy DN-matrix $\mathrm{L}_q^\epsilon$. In the white areas the computation breaks down due to noise and/or large values of $\lambda$. The black line indicates the truncation radius $R_2$ used in $\vct{t}_R$. 

We choose $R_1 = 1.37$ to take care of problems with $\abs{\lambda}$ close to 1. We have $N_\lambda = 256$ as the grid parameter of section \ref{sec:truncation}. We solve the periodized integral equation \eqref{IEIS-modified} in a $2^{M_d}\times 2^{M_d}$ -sized $\lambda$-grid with $M_d=7$.

\subsubsection{Effect of truncation}
\begin{figure}[h]
\begin{picture}(120,150)
\epsfxsize=9cm
\put(17,3){\epsffile{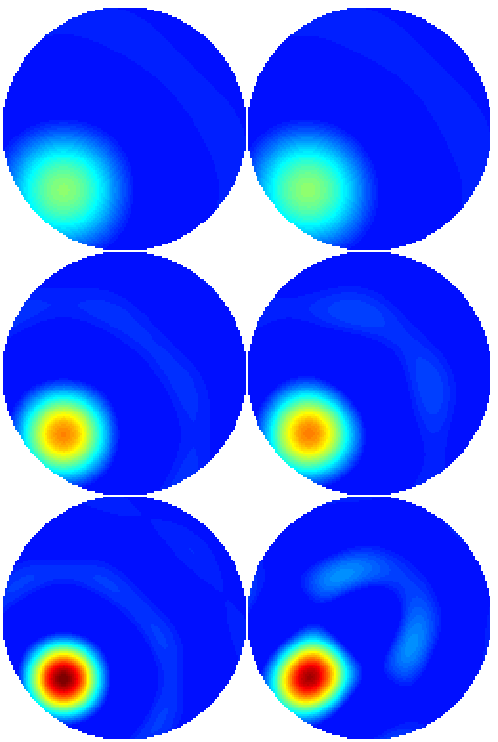}}
\put(0,130){$R_2 = 200$}
\put(0,85){$R_2 = 260$}
\put(0,40){$R_2 = 320$}
\put(26,140){Non-noisy, $\mathrm{L}_q$}
\put(80,140){Noisy, $\mathrm{L}_q^\epsilon$}
\put(22,93){79\%}
\put(67,93){79\%}
\put(22,48){68\%}
\put(67,48){69\%}
\put(22,3){57\%}
\put(67,3){63\%}
\end{picture}
\caption{\label{Rtest} Case 1 reconstructions using three different truncation radii, non-noisy reconstructions on the left and noisy reconstructions on the right. Attached are the relative $L^2$ errors compared to the original potential. See figure \ref{case1recon} for the original potential. Energy level $E=0.001$.}
\end{figure}

Using the Case 1 potential we test different circular truncation radii $R_2$ for the scattering transform $\vct{t}_R$. The result is in figure \eqref{Rtest}, where we see how the reconstruction improves by using a larger truncation radius. 

\subsubsection{Reconstructions}
\begin{figure}[h]
\begin{picture}(120,55)
\epsfxsize=12cm
\put(0,0){\epsffile{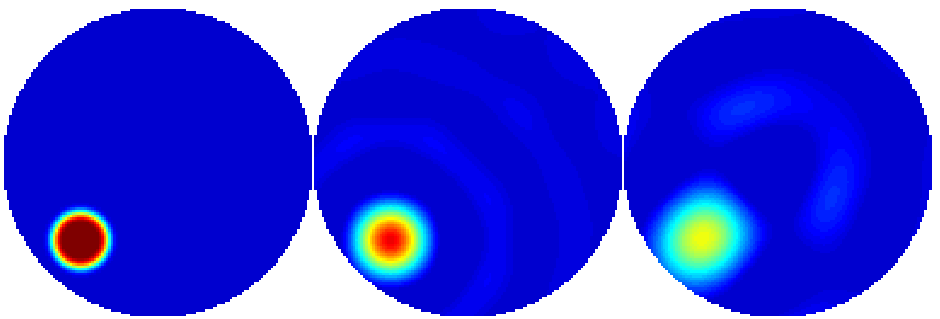}}
\put(6,42){Case 1 potential}
\put(45,42){Reconstruction, $\mathrm{L}_q$}
\put(85,42){Reconstruction, $\mathrm{L}^{\epsilon}_q$}
\put(42,0){40\%}
\put(82,0){63\%}
\end{picture}
\caption{\label{case1recon}On the left: the original Case 1 potential, see figure \ref{fig:q1plot}. In the middle: reconstruction using the non-noisy DN-matrix $\mathrm{L}_q$. On the right: reconstruction using the noisy DN-matrix $\mathrm{L}^{\epsilon}_q$. Relative errors $\norm{q_0-q_{\textrm{rec}}}_{L^2(\Omega)}/\norm{q_0}_{L^2(\Omega)}$ are given. Energy level $E=0.001$. The colormap is different from the one used in figure \ref{Rtest}.}
\end{figure}

\begin{figure}[h]
\begin{picture}(120,55)
\epsfxsize=12cm
\put(0,0){\epsffile{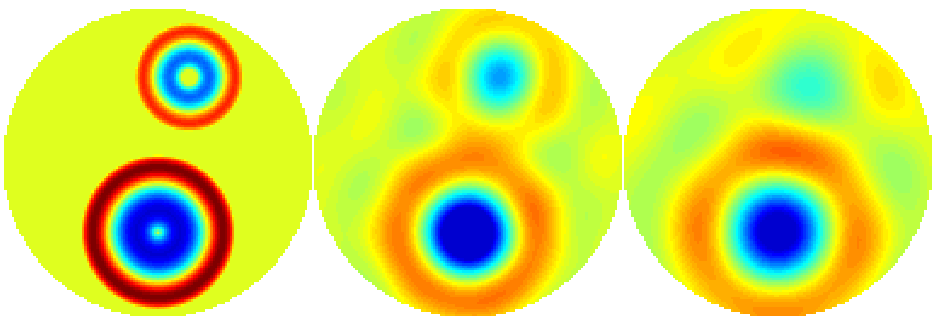}}
\put(6,42){Case 3 potential}
\put(45,42){Reconstruction, $\mathrm{L}_q$}
\put(85,42){Reconstruction, $\mathrm{L}^{\epsilon}_q$}
\put(42,0){61\%}
\put(82,0){73\%}
\end{picture}
\caption{\label{case3recon}On the left: the original Case 2 potential, see figure \ref{fig:q3plot}. In the middle: reconstruction using the non-noisy DN-matrix $\mathrm{L}_q$. On the right: reconstruction using the noisy DN-matrix $\mathrm{L}^{\epsilon}_q$. Relative errors $\norm{q_0-q_{\textrm{rec}}}_{L^2(\Omega)}/\norm{q_0}_{L^2(\Omega)}$ are given. Energy level $E=0.001$.}
\end{figure}

Based on figure \ref{scats2D_DN} and similarly for the other cases, we choose for $R_2=410$ and $R_2=320$ for Case 1 non-noisy and noisy reconstructions from the DN-matrix. For Case 2 we choose $R_2=435$ and $R_2=320$ respectively. In figures \ref{case1recon} and \ref{case3recon} we picture the original potentials on the left, the reconstruction using \eqref{q0recon-discrete} without noise in the middle and the reconstruction using \eqref{q0recon-discrete} with added noise on the right. Relative errors 
$$
\norm{q_0-q_{\textrm{rec}}}_{L^2(\Omega)}/\norm{q_0}_{L^2(\Omega)},
$$ 
where $q_{\textrm{rec}}$ is the reconstruction, are given.

\subsection{Comparison of algorithms}\label{sec:comparison}
\begin{figure}[h]
\begin{picture}(120,200)
\epsfysize=19cm
\put(32,2){\epsffile{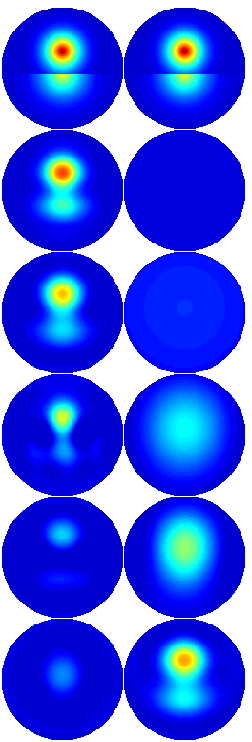}}

\put(10,174){Original $q_0$}
\put(13,142){\fbox{$E=0.1$}}
\put(13,110){\fbox{$E=1.0$}}
\put(13,79){\fbox{$E=5.0$}}
\put(13,47){\fbox{$E=10.0$}}
\put(13,16){\fbox{$E=30.0$}}

\put(23,133){\small $R_2=45$}
\put(23,101){\small $R_2=12$}
\put(23,69){\small $R_2=6.1$}
\put(23,37){\small $R_2=3.7$}
\put(23,5){\small $R_2=1.6$}

\put(37,192){Our method}
\put(64,192){Novikov-Santacesaria}

\put(30,128){15\%}
\put(30,96){27\%}
\put(30,65){59\%}
\put(30,33){88\%}
\put(30,1){89\%}

\put(63,128){98\%}
\put(63,96){83\%}
\put(63,65){53\%}
\put(63,33){39\%}
\put(63,1){20\%}
\end{picture}
\caption{\label{fig:comparison}Numerical comparison of the two algorithms. Top row: the original potential $q_0$. The next rows show the reconstruction using increasing energy $E$, our method on the left, the method of Novikov-Santacesaria on the right. The truncation radii $R_2$ for out method and relative $L^2$ errors are also shown.}
\end{figure}

We used another test potential for the comparison of our algorithm against the Novikov-Santacesaria algorithm \cite{Novikov01012013}. The result is pictured in \ref{fig:comparison}: on top we have the original potential having values between 0 and 1, below we have the reconstructions of both algorithms using non-noisy DN-maps for energies $E=0.1,1,5,10$ and $E=30$. Our method is on the left, the other method on the right. Relative errors and truncation radii for our method are given.

\begin{figure}[h]
\begin{picture}(120,55)
\epsfxsize=12cm
\put(0,5){\epsffile{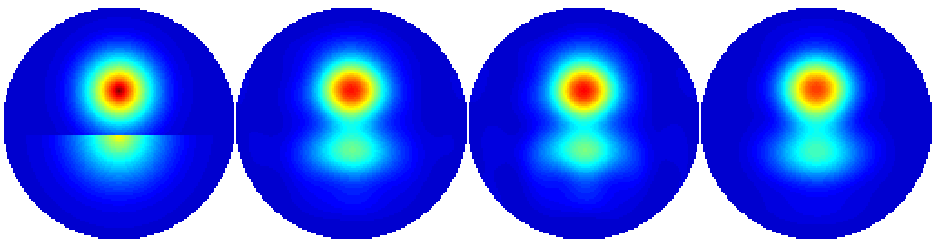}}
\put(6,40){Original $q_0$}
\put(38,42){$E=10^{-5}$}
\put(68,42){$E=10^{-3}$}
\put(98,42){$E=10^{-1}$}
\put(37,38){\small $R_2 = 5000$}
\put(67,38){\small $R_2 = 520$}
\put(97,38){\small $R_2 = 45$}
\put(39,0){12.57\%}
\put(69,0){12.45\%}
\put(99,0){15.13\%}
\end{picture}
\caption{\label{fig:Etest}Reconstruction using three different, small energy levels. Relative errors $\norm{q_0-q_{\textrm{rec}}}_{L^2(\Omega)}/\norm{q_0}_{L^2(\Omega)}$ and truncation radii are given.}
\end{figure}

For smaller energies than $E=0.1$, we get approximately the same reconstruction. In figure \ref{fig:Etest} we use energies $10^{-5}, 10^{-3}$ and $0.1$. The small differences can be attributed to differences in the truncation that we choose based on the plotted scattering transform as in figure \ref{scats2D_DN}.

\section{Conclusions}
\noindent
We developed a new numerical method for reconstructing the potential from boundary measurements in the Gel'fand-Calder\'on problem. The method seems to work as evidenced by the reconstructions, even if the theory is still missing details: the operator $\MM$ is not used in the reconstructions. See the reconstructions of Cases 1 and 2, pictured in figures \ref{case1recon} and \ref{case3recon}. Also see the comparison reconstructions of figures \ref{fig:comparison} and \ref{fig:Etest}. The reconstructions of the potential in the comparisons has significantly lower numerical relative error, which we attribute to the lower contrast of the potential. We conclude that our method works better for lower contrast potentials, even if the reconstructions of Case 1 and 2 are visually satisfactory and do reveal important features of the original potential.

The numerical method for $g_\lambda$ is not accurate near $\abs{\lambda}=1$ which resulted in high errors in the verification test of figure \ref{fig:dbar-error} using the $\dbar$ -equation. For other values of $\lambda$ this numerical method is accurate enough for good quality reconstructions.

Regarding the radial potentials, the numerical evidence show no exceptional points for small $\alpha$ nor for large $\lambda$ which is to be expected according to the theory for small potentials. Also according to our tests there are no exceptional points for positive $\alpha$. For negative $\alpha$, there are either one or two exceptional circles in the range of parameters investigated. The two types of potentials $q_\alpha^{(1)}$ and $q_\alpha^{(2)}$ have little difference in their exceptional points, mainly in the second exceptional circle forming at $\alpha=-20$ as $\alpha$ is decreased from zero.

The comparison result of figure \ref{fig:comparison} shows that currently our method works better with smaller energies and becomes uneffective at larger energies as is expected from the stability results. Increasing the energy scales the usable area of the truncation down, as indicated by the truncation radii used in the comparison. Using energies $10^{-5}\leq E\leq 0.1$ result to virtually the same reconstruction, and as of now this is the optimal range for our method - further fine-tuning is left for future works.

\section*{Acknowledgements}
JPT was supported in part by the Finnish Cultural Foundation and European Research Council (ERC). ML and SS were supported by the Finnish Centre of Excellence in Inverse Problems Research 2012--2017 (Academy of Finland CoE-project 250215). MS was supported by FiDiPro project of Academy of Finland, number 263235.

\bibliographystyle{plain}

\end{document}